\documentclass[11pt,a4paper]{article}

\usepackage[numbers]{natbib}
\usepackage[english]{babel}
\usepackage[top=2in,bottom=1.5in, outer=1.25in,inner=1.5in]{geometry}
\usepackage[bookmarks,breaklinks]{hyperref}

\usepackage[matrix, arrow,all,cmtip,color]{xy}
\usepackage{enumitem}

\usepackage{setspace}

\usepackage[centertags]{amsmath} 
\usepackage{amsthm}
\usepackage{amssymb}

 \usepackage[textsize=tiny]{todonotes}

\theoremstyle{plain} 
\newtheorem{conjecture}{Conjecture}

\newtheorem{thm}{Theorem}[section]
\newtheorem*{thm*}{Theorem}
\newtheorem{cor}[thm]{Corollary}
\newtheorem{fact}[thm]{Fact} 
\newtheorem{lemma}[thm]{Lemma}

\newtheorem{prop}[thm]{Proposition} 
\newtheorem{defn}[thm]{Definition}
\theoremstyle{remark}
\newtheorem*{rem}{Remark}

\newcommand{\A}{\mathbb{A}}

\newcommand{\Z}{\mathbb{Z}}

\newcommand{\G}{\mathbb{G}}

\newcommand{\Ii}{\mathcal{I}} 
 
\newcommand{\Mm}{\mathcal{M}}
\newcommand{\M}{\mathfrak M}
\newcommand{\Oo}{\mathcal{O}} 
\newcommand{\m}{\mathfrak{m}}

\newcommand{\Spec}{\operatorname{Spec}}

\newcommand{\red}{\operatorname{red}}

\newcommand{\acl}{\operatorname{acl}}

\newcommand{\Fr}{\operatorname{Fr}}

\newcommand{\tp}{\operatorname{tp}}
\newcommand{\Cb}{\operatorname{Cb}}

\newcommand{\eps}{\operatorname{\varepsilon}}

\newcommand{\Ker}{\operatorname{Ker}}

\newcommand{\End}{\operatorname{End}}
\newcommand{\Aut}{\operatorname{Aut}}

\newcommand{\what}{\widehat}
 
\newcommand{\wo}{\setminus} \newcommand{\set}[1]{\{\ #1\ \}}
 \newcommand{\tuple}[1]{(\ #1\
)} \newcommand{\suchthat}[2]{\{\ #1\ \mid\ #2\ \}}

\newcommand{\0}{\emptyset}
\def\nsim{\raise.17ex\hbox{$\scriptstyle\sim$}}

\DeclareFontFamily{U}{mathx}{\hyphenchar\font45}
\DeclareFontShape{U}{mathx}{m}{n}{<-> mathx10}{}
\DeclareSymbolFont{mathx}{U}{mathx}{m}{n}
\DeclareMathAccent{\widebar}{0}{mathx}{"73}

\def\ind{\mathrel{\raise0.2ex\hbox{\ooalign{\hidewidth$\vert$\hidewidth\cr\raise-0.9ex\hbox{$\smile$}}}}}

\newcommand{\groupconf}[6]{
   \xy <1.5cm,0cm>: (1,0.2)*{#4}; (2,0.2)*{#5}; (0.8,-0.8)*{#6};
  (-0.2,0)*{#1}; (-0.2,-1)*{#2}; (-0.2,-2)*{#3}; (1,0)*{}; (0,-2)*{}
  **\dir{-}; (0,-2)*{}; (0,0)*{} **\dir{-}; (0,0)*{}; (2,0)*{}
  **\dir{-}; (2,0)*{}; (0,-1)*{} **\dir{-};
  \endxy }

\begin{document}

\bibliographystyle{abbrvnat}

\begin{center}
 {\LARGE\bf Incidence systems on Cartesian powers of algebraic
   curves}\\[2ex] 
\end{center}

\begin{center}
  {\normalsize Assaf Hasson\footnote{was supported by ISF grants
      No.~181/16 and 1156/10.}  and Dmitry Sustretov \footnote{has
      received funding from the European Union's Horizon 2020 research
      and innovation programme under the Marie Sklodowska-Curie grant
      agreement No.~843100 (NALIMDIF) and was supported by the
      European Research Council under the European Community's Seventh
      Framework Programme (FP7/2007-2013) with ERC Grant Agreement
      nr.~615722 (MOTMELSUM).} }
\end{center}
\begin{center}
\begin{minipage}[t]{0.9\linewidth}
\begin{abstract}
  We show that a non locally modular reduct of the full Zariski
  structure of an algebraic curve interprets an infinite field.

\end{abstract}
\end{minipage}
\end{center}

\begin{center}
  \small
  \tableofcontents
\end{center}


\section{Introduction}
In \cite[\S 2]{ArtinGeometry} Artin describes the basic problem of
classifying abstract plane geometries (viewed as incidence systems of
points and lines) as follows ``Given a plane geometry [...] assume that
certain axioms of geometric nature are true [...] is it possible to find
a field $k$ such that the points of our geometry can be described by
coordinates from $k$ and lines by linear equations?''. Zilber's
trichotomy principle (to be described in more detail in the next
section) can be viewed as an abstraction of the above problem,
replacing the ``axioms of geometric nature'' with a well behaved
theory of dimension (see, e.g., \cite[\S 1]{poizat}).

Conjectured in various forms by Zilber throughout the late 1970s,
essentially every aspect of Zilber's trichotomy, in its full
generality, was refuted by Hrushovski \cite{Hr1}, \cite{Hr2} in the
late 1980s. Due to Hrushovski's cornucopia of counterexamples the
conjecture has never been reformulated. Yet, Zilber's principle
remains a central and powerful theme in model theory: it has been
proved to hold in many natural examples such as differentially closed
fields of characteristic $0$, algebraically closed fields with a
generic automorphism, o-minimal theories and more (see \cite{chatzhr,
  pillay-ziegler, peterzil1998trichotomy, chatzidakis2002model,
  hz}). Many of these special cases of Zilber's trichotomy had
striking applications in algebra and geometry (see
\cite{hrushovski1996mordell,hrushovski2001manin,scanlon2006local}). More
recently, in \cite{ZilSchanuel}, Zilber outlines a model theoretic
framework for studying far reaching extensions of the Mordell-Lang
conjecture. One of the key features of Zilber's strategy is the
trichotomy theorem for Zariski Geometries \cite{hz}.

The key to the classification of Desarguesian plane geometries (the
fundamental theorem of projective geometry) is the reconstruction of
the underlying field $k$ as the ring of direction preserving
endomorphisms of the group of translations.  The reconstruction of a
field out of abstract geometric data is also the essence of Zilber's
trichotomy and is the engine in many of its applications.  A
relatively recent application of one such result is Zilber's model
theoretic proof \cite{ztorelli} of a significant strengthening of a
theorem of Bogomolov, Korotiaev, and Tschinkel \cite{jacobian}. The
model theoretic heart of Zilber's proof is Rabinovich's trichotomy
theorem for \emph{reducts} of algebraically closed fields
\cite{rabinovich}. In the concluding paragraph of the introduction to
\cite{ztorelli} Zilber writes: ``It is therefore natural to aim for a
new proof of Rabinovich' theorem, or even a full proof of the
Restricted Trichotomy along the lines of the classification theorem of
Hrushovski and Zilber \cite{hz}, or by other modern methods
[...]. This is a challenge for the model-theoretic community.'' \\

The conjecture referred to in Zilber's text above can be formulated as
follows\footnote{The content of Conjecture \ref{RTC} is explained for
  non-experts in Section \ref{background}.}$^,$\footnote{In the full
  conjecture referred to by Zilber $\M$ is the full Zariski structure
  on an $n$-dimensional constructible set (for $n$ possibly greater
  than $1$).}:

\begin{conjecture}\label{RTC}
  Let $\Mm$ be a strongly minimal reduct of the full Zariski structure
  on an algebraic curve $M$ over an algebraically closed field $K$
  which is not locally modular. Then there exist $\Mm$-definable
  $L,E$ such that $E\subseteq L\times L$ is an equivalence relation
  with finite classes and $L/E$ with the $\Mm$-induced structure is a
  field $K$-definably isomorphic to $K$.
\end{conjecture}

Rabinovich \cite{rabinovich} proved Conjecture \ref{RTC} in the
special case where $M=\A^1$, and her result can be extended by general
principles to any rational curve. In the present paper we prove
Conjecture~\ref{RTC}. Our approach to the problem follows a, by now,
well known strategy introduced by Zilber and Rabinovich and owes to
\cite{marker-pillay}. Using a standard model theoretic technique,
Hrushovski's field configuration (see Section \ref{S:gc} for details),
the problem is reduced to showing that tangency is (up to a finite
correction) reduct-definable in families.

To achieve this goal, we proceed in two steps. In the first step
(carried out in Section~\ref{tangency}) we study slopes of families of
branches (at a given point) and their behaviour under composition of
curves and, in case the ambient structure is an expansion of a group,
under point-wise addition. This culminates in
Proposition~\ref{tangency-intersection}, which is the key to the
definability of tangency, and in Lemma~\ref{many-slopes}, providing us
with the (algebraic) group which is the template allowing us to
construct the group configuration.

Section~\ref{interp-field}, where Conjecture~\ref{RTC} is proved
(Theorem~\ref{main-theorem}), is dedicated, mainly, to verifying that
the assumptions of the technical result of the previous section can be
met in the reduct. In Section~\ref{sect:unram} we show that our
definition of slope is meaningful in reduct-definable families of
curves (in positive characteristic). In Section~\ref{int-group}, where
the main step towards proving Theorem~\ref{main-theorem} is carried
out, the key difficulty to overcome is in the application of
Proposition~\ref{tangency-intersection}.

As already mentioned, the general scheme of our proof seems to have
much in common with Rabinovich's original work, though we were unable
to understand significant parts of her argument, which are highly
technical. For that reason we cannot pinpoint the reason for the
present work being more general, considerably shorter, and technically
simpler.

Finally, it should be mentioned that the tools developed in the
present paper seem to extend naturally to various other contexts. For
example, one can envisage extending the results of \cite{piotr-serge}
to positive characteristic, and any algebraic group and -- possibly --
even a full proof of the restricted trichotomy conjecture for
structures definable in ACVF (at least modulo the problem of showing
that the 1-dimensional group reconstructed by our methods embeds in an
algebraic group: it can probably be shown that the group will always
be isomorphic to either $\mathbb G_a$ or to $\mathbb G_m$).

\section{Model-theoretic background}
\label{background}
For readers unfamiliar with model theory we give a self contained
exposition of Conjecture \ref{RTC}. In order to keep this introduction
as short as possible, we specialise our definitions to the setting in
which they will be applied. We refer interested readers to
\cite[\S1.1-2]{vdDries} for a more detailed discussion of structures
and definable sets. Readers familiar with the basics of model theory
are advised to skip the remainder of the present section.

Given an algebaraic curve $M$ over an algebraically closed field $k$
(reduced, but not necessarily irreducible, smooth or projective), the
\emph{full Zariski structure on $M$}, denoted $\M$, is the set of
$k$-rational points, $M(k)$ equipped with the collection of all
Boolean algebras of constructible sets on the Cartesian powers
$M^n(k)$. The full Zariski structure on a curve $M$ is an example of
the model theoretic notion of a structure.

A \emph{first-order structure} or simply a \emph{structure}
$\mathcal N$ is a non-empty set $N$ (called the universe of
$\mathcal N$) equipped with a collection, $\mathrm{Def}(\mathcal N)$,
of Boolean algebras
$\mathrm{Def}_l(\mathcal N)\subseteq \mathcal P(N^l)$ for all $l>0$,
such that $\mathrm{Def}_l(\mathcal N)$ contains all diagonals
$\Delta_{i,j}^l:=\{(x_1,\dots x_l): x_i=x_j\}$, and such that
$\mathrm{Def}(\mathcal N)$ is closed under finite cartesian products
and projections of the form
$(x_1,\dots, x_n)\mapsto (x_1,\dots x_{n-1})$.  Somewhat analogously
to geometric terminology the tuples
$(x_1, \ldots, x_n) \in S \subset M^l$ are called \emph{points of the
  definable set $S$}.  If $A\subseteq N$ is any set, a subset
$X\subseteq N^l$ is \emph{definable with parameters in $A$} (or
\emph{$A$-definable}) if there exists a definable set
$Y\subseteq N^{n+m}$ (some $m\ge 0$) such that
$Y=Y_a:=\{x\in N^l: (x,a)\in Y\}$ for some $a\subseteq A$.

Note that by Chevalley's theorem (see, e.g., \cite[Corollary
3.2.8]{marker}), over an algebraically closed field $k$, the
collection of constructible sets on cartesian powers of an algebraic
curve $M$ is closed under projections, and therefore the full Zariski
structure, $\M$, on $M$ is, indeed, a structure in the above sense.
It is a well known fact (e.g., \cite{hz}) that the field $k$ can be
reconstructed from $\M$. Let us now explain more precisely what is
meant by that.

A (partial) function $f:N^l\to N$ is definable if its graph is. Thus,
for example, we say that a group is definable in $\mathcal N$, if
there exists a definable set $G\subseteq N^l$ and a definable function
$p:G\times G\to G$ such that $(G,p)$ is a group (note that the
function $x\mapsto x^{-1}$ is automatically definable if $(G,p)$ is a
group). The definability of a field in a structure $\mathcal N$ is
defined analogously. It is not hard to check (and follows from the
main result of \cite{hz}) that if $\M$ is the full Zariski structure
on an algebraic curve $M$ over an algebraically closed field $k$ then
a field $F$ is definable in $\M$ (and $F$ is isomorphic, definably in
the standard field structure on $k$, to $k$).

But we need a somewhat subtler notion than definability. Consider, as
a simple example, the structure $\mathcal C$ with universe
$\mathbb C\times \{0,1\}$, and whose definable sets are all those of
the form $\{((x_1,i_1), \dots, (x_n,i_n): (x_1,\dots, x_n)\in D\}$
where $D$ is a constructible subset of $\mathbb C^n$ and
$i_j\in \{0,1\}$ for all $1\le j\le n$. It is easy to verify that all
functions definable in $\mathcal C$ are locally constant, and
therefore there is no definable field in $\mathcal C$. Consider,
however, the equivalence relation $x\sim y$ (in $\mathcal C$) defined
by $y\in (1,0)\cdot x$ (recalling the interpretation of multiplication
in $\mathcal C$, this is a $\mathcal C$-definable way of saying that
$x$ and $y$ have the same first coordinate). Then $\sim$ is a
$\mathcal C$-definable equivalence relation, and $\mathcal C/\sim$ is
naturally isomorphic to the full Zariski structure on $\mathbb C$.

In model-theoretic terms the structure $\mathcal C$ in the previous
example \emph{interprets} a field definably isomorphic to
$\mathbb C$. In general, if $\mathcal N$ is a structure, $E$ a
definable equivalence relation on $N^l$ and $\pi: N^l\to N^l/E$ is the
natural projection, the induced structure on $N^l/E$ is the
push-forward of the Boolean algebras on powers of $N^l$ via $\pi$. We
say that $\mathcal N$ \emph{interprets} a field if a field is
definable in the structure induced on $N^l/E$ for some $l$ and
$\mathcal N$-definable equivalence relation $E$ on $N^l$.

In the above example the universe $\mathbb C\times \{0,1\}$ of
$\mathcal C$ is definable in the full Zariski structure on
$\mathbb C$, and every definable set in $\mathcal C$ is definable in
$\mathbb C$. But, as we have seen, $\mathcal C$ is not the full
Zariski structure on $\mathbb C$. The structure $\mathcal C$ is an
example of a \emph{reduct} of the full Zariski structure on
$\mathbb C$. Generally, if $\mathcal M$ is a structure whose universe
is an algebraic curve $M$ and every $\mathcal M$-definable set is
$\M$-definable then $\Mm$ is a \emph{reduct} of $\M$.

Zilber's conjecture is concerned with the question of interpreting a
field in a reduct, $\Mm$, of the full Zariski structure, $\M$, on an
algebraic curve $M$. Assume that an infinite field $\mathbb F$ is
interpretable in $\Mm$. Then by \cite[Theorem 3.2.20]{marker} the
universe $F$ of $\mathbb F$ can be identified with a constructible
subset of $k^l$ for some $l$, and by \cite[Theorem 4.15]{poizat} $k$ is
definably isomorphic to $\mathbb F$. Thus there is a definable
finite-to-finite correspondence $\Psi\subseteq F\times M$. It is easy
to check that $\Psi$ can be taken to be $\Mm$-interpretable (e.g., if
$F$ is definable in $\Mm$ then $\Psi$ can be taken to be the graph of
a projection function, the general case is slightly more delicate and
we skip the details). If we push the family of affine lines in $F^2$
via $\Psi$ we obtain a 1-dimensional constructible subset $U$ of $M$
such that for any $p,q\in U$ there is a curve $C:=\Psi(L)$ -- for $L$
an affine line in $F^2$ -- with $p,q\in C$. We have thus verified that
for $\Mm$ to interpret a field it is necessary that there exists a
2-dimensional constructible $U\subseteq M^2$ and a definable set
$X\subseteq M^{2+l}$ such that $X_t:=\{(x,y): (x,y,t)\in X\}$ is
1-dimensional (or empty) for all $t\in M^l$ and such that for all
$p,q\in U$ there exists $t\in M^l$ such that $p,q\in X_t$. The main
result of the present work, Theorem~\ref{main-theorem}, states that
this condition is, in fact sufficient.

\begin{defn}\label{ample}
  Let $\Mm$ be a reduct of the full Zariski structure $\M$ on an
  algebraic curve $M$ over an algebraically closed field $k$. An
  $\Mm$-definable \emph{ample family of curves in $M^2$} is a set
  $X\subseteq M^{2+l}$ such that
  \begin{itemize}
  \item $\dim(X_t)=1$ for all $t\in M^l$ such that $X_t\neq \0$ and
  \item there exists a 2-dimensional $U\subseteq M^2$ such that for
    all $p,q\in U$ there exists $t\in M^l$ with $p,q\in X_t$.
  \end{itemize}
\end{defn}

In model-theoretic terms the existence of an ample family as above is
equivalent, \cite[Lemma 8.1.13]{marker},  to \textbf{\emph{non local
    modularity}} of the structure $\mathcal M$. \\ 

If $X$ is an ample family in $M^2$ we denote by $(M,X)$ the smallest reduct of $\M$ containing $X$. We can thus formulate Conjecture~\ref{RTC}: 
\begin{conjecture}[Zilber's restricted trichotomy in dimension
	1]\label{conjB}
	Let $M$ be an algebraic curve over an algebraically closed field
	$k$. Let $X\subseteq M^2\times T$ be the total space of an ample
	family in $M^2$,
	then a field $K$ is interpretable in $\mathcal M=(M,X)$.
\end{conjecture}

In \cite[\S 2.4]{ArtinGeometry} not only is the field recovered from
the affine geometry, but also the geometry is recovered as the affine
plane over that field. In the present setting, there are examples due
to Hrushovski (see, e.g., \cite{Martin}) showing that the full Zariski
structure of the curve $M$ cannot be recovered from $\Mm$. This can
probably be achieved if $X$ is \emph{very ample} in the sense of
\cite{hz} (namely, if the set $X$ in Definition \ref{ample} the separates points in $M^2$), but we do not
study this question here.

\section{Tangency}
\label{tangency}
The reconstruction of the field is obtained in two steps. First, we
reconstruct a 1-dimensional algebraic group, and then -- using the
group structure to sharpen the same arguments -- we reconstruct the
field. Roughly, the reconstruction of a group is obtained in three
stages: first we identify a reduct definable family $X\to T$ of
algebraic curves whose associated family of slopes at some point
$P=(a,a)\in M^2$ is a 1-dimensional algebraic group under
composition. The second, and most crucial part of the proof is
commonly dubbed \emph{definability of tangency.} In its cleanest form
this consists in showing that, given families $X\to T$ and $Y\to S$ as
above, the set of all $(t,s)\in T\times S$ such that $X_t$ is tangent
(in an appropriate sense) to $Y_s$ at $P$ is $\Mm$-definable. Finally,
the group is reconstructed by invoking \emph{the group configuration
  theorem}, a well known model theoretic technique (to be described in
more detail in the next section), using the results of the previous
stages. In the next two subsections we take care of the two first
stages of the this strategy.

Before starting with our set up on the technical level let us discuss
some of the challenges that motivated the definitions to be shortly
presented. In the implementation of the strategy outlined above two
 difficulties arise.

 Firstly, if we consider only the first-order slopes, then due to
 inseparability issues in positive characteristic it becomes hard to
 find a 1-dimensional family of curves definable in the reduct such
 that its associated slopes at some point range in a 1-dimensional set
 --- such a family is needed to construct the first group
 configuration (Section~\ref{int-group}). The solution is to consider
 tangency information up to any order $n$ and pick the order so that
 there are enough slopes. Interestingly --- and this was apparent
 already in \cite{marker-pillay} --- in the presence of a group
 structure, the problem doesn't arise, which is a good coincidence,
 since the second group configuration (Section~\ref{int-field}) has to
 be built using the first-order tangency information.

 Secondly, we can't work only with smooth points to define the slope,
 since the operations of composition and point-wise addition that are
 used in the construction of the group configurations do not preserve
 smoothness. Our approach to this is to track a particular branch of a
 curve at a particular point as the operations of composition and
 point-wise addition are applied: one can then have control over the
 slope of a branch, appropriately defined. Note that we use the term
 \emph{branch} (Definition~\ref{branches}) in a more restrictive sense
 than what is usually understood by it: in a way our branches are
 `always smooth' (or more precisely `always \'etale over the first
 factor $M$ of $M^2$'), so that the notion of a slope always makes
 sense for them. For any curve in $M^2$ the projection either on the
 first or on the second factor $M$ is going to be generically \'etale
 (Lemma~\ref{projections}), even in positive characteristic, and so
 there is going to be a unique branch at any general enough point on
 this curve (Lemma~\ref{etale-branch}). This statement
 generalizes appropriately to families of curves too. By virtue of
 Propositions~\ref{compo-slope} and \ref{sum-slope} the slopes of
 relevant branches can be tracked as the curves are composed and
 point-wise added. All curves and branches that we work with in
 Section~4 are obtained this way.

\subsection{Slopes and operations on correspondences}
\label{slopes}

First of all, our main objects of interest are definable subsets in a
reduct of the full Zariski structure on a fixed curve $M$, we will
adopt the following non-standard terminology: we will call a
constructible subset $Z\subseteq M^n$ (for some $n$) of dimension 1 a
\emph{curve} even if it is reducible and if it has connected
components of dimension 0. If a curve $Z$ does not contain connected
components of dimension 0 then we call $Z$ a \emph{pure-dimensional
  curve}. Clearly, every curve contains a maximal pure-dimensional
curve. In the few situations when we refer to abstract algebraic
curves (that is, algebraic varieties of pure dimension 1 over a fixed
algebraically closed field) we will use the term \emph{algebraic
  curve}. We will not distinguish notationally between subsets of
$M^n$ definable in a reduct of the full Zariski structure on $M$ and
constructible subsets of the varieties (or even schemes) $M^n$, and in
particular between definable curves and their algebro-geometric
counterparts.

Recall that any algebraic variety over a perfect field admits a dense
Zariski open subset that is smooth (see, e.g., Corollary to
Theorem~30.5 of \cite{matsumura1989commutative}). Let $Z \subset M^2$
be a pure-dimensional curve, and $a=(a_1, a_2) \in Z$ be a smooth
point of $M^2$. Since the completion of the local ring of a smooth
point of a variety is a formal power series ring
\cite[Theorem~29.7]{matsumura1989commutative}, we can pick
some isomorphisms
$$
\widehat{\Oo_{M, a_1}} \cong k[[x]] \qquad \widehat{\Oo_{M, a_2}}
\cong k[[y]]
$$
inducing an isomorphism $\widehat{\Oo_{M^2,a}} \cong k[[x,y]]$, and
then $\widehat{\Oo_{Z,a}} = k[[x,y]]/(f)$ for some $f \in k[[x,y]]$.
We call \emph{branches of $Z$ at $a$} the factors in the prime
decomposition of $f$ of the form $y - g_\alpha, g_\alpha \in k[[x]]$
(note that this is different from the standard use of the term, but
since we will never use the term in its standard meaning in this
article, no confusion will occur). In particular, if the projection of
$Z$ onto the first factor $M$ in $M^2$ is \'etale in a neighbourhood
of $a$, by Hensel's lemma (stated as in \cite[\S 4,
Theorem~4.2(d)]{milne}), the natural morphism
$k[[x]] \to k[[x,y]]/(f)$ is an isomorphism and therefore $f$ can be
written uniquely as $u(y - g)$ where
$u \in k[[x]]^\times, g \in k[[x]]$. We call the truncation to the
$n$-th order of the series $g_\alpha$ the \emph{$n$-th order slope of
  a branch $\alpha$ of $Z$}. Naturally, the slope of a branch of a
pure-dimensional curve depends on the choice of the isomorphism
$\widehat{\Oo_{M^2,a}} \cong k[[x,y]]$, but this choice does not
affect any of the properties of slopes we will be interested in.

We view curves in $M^2$ as finite-to-finite correspondences between
the two factors $M$. The purpose of the present section is to study
the behaviour of slopes of branches with respect to two natural
operations on correspondences: composition and ``point-wise addition''
(see Definition~\ref{sum}) when $M$ has a structure of an algebraic
group. We will show that if $Z, W$ are two curves and $\alpha, \beta$
are two branches of $Z,W$ at $a=(a_1,a_2), b=(b_1,b_2) \in M^2$
respectively, and $a_2=b_1$ then the composition $W \circ Z$ has a
branch $\beta \circ \alpha$ at $(a_1,b_2)$ such that its slope is the
composition of the $n$-th order slopes of $\alpha$ and $\beta$ (as
truncated polynomials) whenever the latter are defined. A similar
statement can be made about the slopes of branches of curves that are
`point-wise added' if $M$ has a structure of a group. Later we will
construct a group configuration starting from a family of curves
definable in a reduct of a full Zariski structure on $M$ such that the
set of its $n$-th order slopes at a given point coincides, up to a
finite set, with a 1-dimensional algebraic subgroup of
$\Aut(k[[x]]/(x^{n+1}))$ (a truncated polynomial $f$ corresponds
naturally to the automorphism of $k[[x]]/(x^{n+1})$ sending $x$ to
$f$). 

Since we will have to work with families of curves, we will also
introduce notions of families of branches and slopes. When the
characteristic of the base field is positive, we will often have to
work with curves and families of curves in $M \times M^{(p^n)}$ where
$M^{(p^n)}$ is the pull-back of $M$ by the Frobenius endomorphism on
$k$ (see Section~\ref{sect:unram}). That is why we do not assume that
factors of the ambient product variety are isomorphic in the
definitions below.

\begin{defn}[Families of curves]
  \label{families}
  If $X_1, X_2$ are two algebraic curves then by a \emph{family of
  	pure-dimensional curves in $X_1 \times X_2$} we will understand a
  finite union $Z$ of pure codimension 1 locally closed subsets
  $Z_i \subset X_1 \times X_2 \times T$, where $T$ is a variety, such
  that $Z_t$ is a pure-dimensional curve for all $t \in T$.  By a
  \emph{family of curves} we understand a constructible subset
  $Z \subset X_1 \times X_2 \times T$ where $T$ is a constructible
  subset of a variety and such that $Z_t$ is a curve for all
  $t \in T$. If $X_1 = X_2 = M$ and $T \subset M^l$ for some $l$ and
  $Z$ is definable in a reduct of the full Zariski structure on $M$ we
  call it a \emph{definable family of curves}.
\end{defn}

While families of curves arise naturally in definable context, in
order to apply the machinery of slopes we need to work with families
of pure-dimensional curves. As long as $T$ is a variety, a family of
curves $Z \subset X_1 \times X_2 \times T$ contains a unique maximal
family of pure-dimensional curves.  The total space $Z$ of a family of
pure-dimensional curves is not necessarily a variety; while this is a
desirable property that will be important in
Subsection~\ref{flat-families}, we do not include it in the definition
so that it can be readily seen that the operations of composition and
point-wise addition preserve the class of families of pure-dimensional
curves.  However, one can easily ensure that the total space is a
variety at the cost of shrinking the parameter space.

\begin{lemma}
  \label{pure-dim}
  Let $T$ be a constructible subset of a variety,
  $W \subset X_1 \times X_2 \times T$ be a family of curves. Then
  there exists a Zariski dense subset $T' \subset T$ which is a
  variety and a maximal locally closed $W' \subset W \times_T T'$
  which is a family of pure-dimensional curves. In particular, $W'$ is
  a variety.
\end{lemma}

\begin{proof}
  It is easily checked using Noetherian induction that any
  constructible set contains a Zariski dense locally closed subset. In
  particular, there exists a dense locally closed subset
  $T_0 \subset T$ which is therefore a variety. Without loss of
  generality we may assume $T_0$ connected. Then $W \times_T T_0$ is a
  union of locally closed sets of the form $W_i \wo Z_i, i \in I$
  where $W_i$ and $Z_i \subset W_i$ are Zariski closed and distinct,
  and the index set $I$ is finite.  Let $I' \subset I$ be the set of
  those indices $i$ for which $W_i$ is codimension 1 in
  $M^2 \times T$. Further shrinking $T_0$ we may assume that
  $\dim (W_i)_t=1, \dim (Z_i)_t=0$ for all $t \in T_0$ and all
  $i \in I'$ (in particular, $Z_i \neq \emptyset$). We put
  $W' = \cup_{i \in I'} W_i \wo Z_i$. It now suffices to show that if
  $W_1 \wo Z_1$ and $W_2 \wo Z_2$ are as above then there exists a
  dense open $T' \subset T_0$ such that
  $(W_1 \wo Z_1) \times_T T' \cup (W_2 \wo Z_2) \times_T T'$ is
  locally closed, the statement of the Lemma then follows by induction
  on the size of $I'$.

  We have
  $$
  (W_1 \wo Z_1) \cup (W_2 \wo Z_2) = (W_1 \cup W_2) \wo \left((Z_1
    \cap Z_2) \cup (Z_1 \cap (W_2 \wo Z_2)) \cup (Z_2 \cap (W_1 \wo Z_1))\right)
  $$
  It follows from an easy dimension computation that
  $$
  \dim (Z_1 \cap (W_2 \wo Z_2)) < \dim  Z_1  \qquad
   \dim (Z_2 \cap (W_1 \wo Z_1)) < \dim Z_2
  $$
  so in particular the projections of $Z_1 \cap (W_2 \wo Z_2)$,
  $Z_2 \cap (W_1 \wo Z_1)$ to $T_0$ are not dense. If $T'$ is a
  dense open set in the complement of the projections then
  $$
  (W_1 \wo Z_1) \times_T T' \cup (W_2 \wo Z_2)\times_T T' = (W_1 \cup
  W_2)\times_T T' \wo (Z_1 \cap Z_2)\times_T T'
  $$
  which is locally closed.  
\end{proof}

Given a family of pure-dimensional curves $Z$ as above, we would like
to be able to pick branches of the curves $Z_t$ depending
algebraically on the parameter $t \in T$. In this case the local
equation of $Z$ in a formal neighbourhood of $\{a\} \times T$ may only
exist locally on $T$, and in order to capture this idea we have to
phrase the definition in terms of formal schemes (we refer the reader
to \cite[II.9]{hartshorne} or any other standard algebraic geometry
reference for the definition of formal schemes).

\begin{defn}[Branches and families of branches]
  \label{branches}
  Let $Z \subset V:=X_1 \times X_2 \times T$ be a family of
  pure-dimensional curves, $a \in M^2$ and assume that $a \in Z_t$ for
  all $t \in T$. Let $\hat{X_1}$ be the formal completion of
  $X_1 \times T$ along $\{a_1\} \times T$, and let $\hat{Z}$ be the
  formal completion of $Z$ along $\{a\} \times T$. A \emph{family of
    branches of $Z$ at $a$} is a closed formal subscheme
  $\hat{Z}_\alpha$ such that the natural projection
  $\hat{Z}_\alpha \to \hat{X_1}$ is an isomorphism. We will call local
  generators of the ideal sheaf that defines $\hat{Z}_\alpha$
  \emph{local equations of $\alpha$}. When $Z \subset X_1 \times X_2$
  is a single curve, we regard it as a family parametrized by a single
  point, and we call families of branches of $Z$ just branches. Given
  a family of branches $\alpha$ we will denote $\alpha_t$ the branch
  given by the fibres $\hat{Z}_{\alpha_t}$ for all $t \in T$.
\end{defn}

\begin{rem}
  If $Z \subset X_1 \times X_2 \times T$ is a family of curves and $T$
  is a variety then in order to simplify the exposition we will refer
  to branches of $Z$ meaning branches of a family of pure-dimensional
  curves $Z_0 \subset Z$.
\end{rem}

Let $X_i$ be algebraic curves, and let
$a=(a_1, \ldots, a_n) \in X=X_1 \times \ldots \times X_n$ be a smooth
point. We say that a \emph{local coordinate system at $a$} is picked
when an isomorphism $\widehat{\Oo_{X_i,a_i}} \cong k[[x_i]]$ is picked
for each $a_i$; in this case we understand that there exists an
isomorphism $\Oo_{X,a} \cong k[[x_1, \ldots, x_n]]$ induced by these
isomorphisms. If local coordinate systems are picked at
$a=(a_1, a_2) \in X_1 \times X_2$, $b=(b_1, b_2) \in X_2 \times X_3$,
we understand without explicit mention that local coordinate systems
are automatically picked at the points
$(a_1, b_2) \in X_1 \times X_3, (b_2, a_1) \in X_3 \times X_1$ which
will be of interest to us later on. Similarly, if $X_1$ has a group
structure and a local coordinate system is picked at a point
$a \in X_1$ then we assume it picked at any point $a' \in X_1$ via
translation. The next lemma gives a sufficient condition for the
existence of a family of branches at a point.

Recall that a morphism of schemes $f: X \to Y$ is called
\emph{quasi-finite} if the fibres $f^{-1}(y)$ are finite for all
$y \in Y$. A quasi-finite morphism $f: X \to Y$ of locally Noetherian
schemes is \emph{unramified} if $\Omega_{X/Y}=0$ (see \cite[Ch.~6,
Corollary~2.3]{liu}) where $\Omega_{X/Y}$ is the module of K\"ahler
differentials of the morphism $f$. A morphism $f$ locally of finite
type is called \emph{\'etale} if it is flat and unramified. Basic
properties of these notions will be recalled in detail and with
references in Section~\ref{flat-families}.

\begin{lemma}
  \label{etale-branch}
  If $Z \subset V=X_1 \times X_2 \times T$ is a family of
  pure-dimensional curves and the projection $Z \to X_1$ is \'etale in
  a neighbourhood of $\{a\} \times T$ for some $a \in X_1 \times X_2$,
  then there exists a unique family of branches of $Z$ at $a$.
\end{lemma}

\begin{proof}
  For any affine open $\Spec R \subset X_1 \times T$ let
  $\Spec R' \subset Z$ be an affine open \'etale over $\Spec R$, let
  $I, I'$ be the ideals vanishing on $\{a_1\} \times T$,
  $\{a\} \times T$ respectively, $\what{R}$ and $\what{R'}$ their
  respective completions. Then by
  \cite[\href{https://stacks.math.columbia.edu/tag/0ALJ}{Tag
    0ALJ}]{stacks} $(\what{R},I)$ is a Henselian pair and by
  \cite[\href{https://stacks.math.columbia.edu/tag/09XI}{Tag
    09XI}]{stacks} there exists a unique isomorphism
  $\what{R'} \to \what{R}$ that defines the unique family of branches.
\end{proof}

\begin{defn}[Slope]
  \label{relslope-definition}
  Let $X_1, X_2$ be algebraic curves, $Z \subset V:=X_1 \times X_2$ a
  pure-dimensional curve, $a \in V$ a smooth point, $a \in Z$, and
  $I, \m_a \subset \Oo_{V,a}$ the ideals of functions that vanish on
  $Z, \{a\}$, respectively.  Assume that a local coordinate system is
  chosen at $a$, so that
  $\varprojlim \Oo_{V,a}/\m_a^n \cong k[[x,y]]$. A branch $\alpha$ of
  $Z$ is therefore defined by a principal ideal $J$ with the property
  that the composition $k[[x]] \to k[[x,y]] \to k[[x,y]]/J$ is an
  isomorphism. The inverse of this isomorphism sends $y$ to
  $f \in x\,k[[x]]$, and $y - f \in J$. We call
  $f \mod x^{n+1} \in k[x]/(x^{n+1})$ the \emph{$n$-th order slope of
    $Z$ at $\alpha$}, denoted $\tau_n(Z, \alpha)$.
\end{defn}

Note that $\tau_n(Z, \alpha)$ depends on the choice of the local 
coordinate system at $a$, and that if an $n$-th order slope of $Z$ at 
$\alpha$ is defined, then the slopes of all orders of $Z$ at $\alpha$ 
are defined.

\vspace{1ex}
\begin{rem} \,\\
  \noindent\begin{enumerate}[label=(\roman*)]
  \item Let $f, g$ be local equations of branches $\alpha, \beta$ at a
    point $a$ of pure-dimensional curves $Z_1, Z_2$, respectively.  If
    $\tau_n(Z_1, \alpha) = \tau_n(Z_2, \beta)$ but
    $\tau_{n+1}(Z_1, \alpha) \neq \tau_{n+1}(Z_2, \beta)$, then
    $f \equiv g \mod x^{n+1}$ and therefore $f - g = x^{n+1}\cdot r$
    for some unit $r \in k[[x]]$ and so
  $$
  k[[x,y]]/(y - f, y - g) \cong k[[x]]/(x^{n+1}).
  $$
  In particular, if $Z_1, Z_2$ are smooth at $a$ and $\alpha, \beta$
  are their unique respective branches at $a$, then the intersection
  multiplicity of $Z_1$ and $Z_2$ at $a$ (as defined in, for example,
  \cite[ex.~I.5.4]{hartshorne}) is $n$.
\item Let $X \subset M^2 \times T$, $Y \subset M^2 \times S$ be
  families of pure-dimensional curves such that $a \in X_t \cap Y_s$
  for all $t \in T, s \in S$ and $X_t \cap Y_s$ is zero-dimensional for
  generic $t,s$.  Let $\alpha, \beta$ be some families of branches of
  $X$ and $Y$ at $a$. Then it follows from Krull's maximal ideal
  theorem that there exists a maximal integer $n$ such that
  $$
  \tau_n(X_t,\alpha_t) = \tau_n(Y_s, \beta_s)
  $$
  for all $t \in T, s \in S$.
\end{enumerate}
\end{rem}

For the benefit of the reader we explain what data in the
Definitions~\ref{branches} and \ref{relslope-definition}
specifies families of branches and slopes, specialising the
description in the language of formal schemes to an affine situation.
Take Zariski open subsets $U \subset X_1 \times X_2$, $W \subset T$ such
that $a \in U$. Let $S, R$ be the rings of regular functions on $U, W$
and let $J_a,J \subset R \otimes S$ be the ideals of regular functions
that vanish on $\{a\} \times W$, $Z \cap U \times W$, respectively. We
fix a local coordinate system at $a$ which gives an isomorphism
$\varprojlim (R \otimes S)/J_a^n \cong R[[x,y]]$.  A choice of a
family of branches $\alpha$ is a choice of an element
$f_\alpha \in R[[x]]$ such that $y - f_\alpha \in R[[x,y]]$ generates
an ideal (necessarily prime) containing $J R[[x,y]]$. The slope
$\tau_n(Z_t, \alpha_t)$ is the truncated polynomial
$f_\alpha \otimes k(t) \mod x^{n+1}\in k[[x]]/(x^{n+1})$. From this
description it is clear that if we regard the $n$-th order slope of
$Z$ at $\alpha_t$ as a tuple of coefficients of
$f_\alpha \otimes k(t)$, then $t \mapsto \tau_n(Z_t, \alpha_t)$ is a
regular function from $W$ to $\A^n$.

Note that the notion of slope is invariant under extensions of the
base field. Assume that all objects in the previous paragraph are
defined over $k$, and let $k' \supset k$ be a field extension. Then
there exists a family of branches $\alpha_{k'}$ of
$Z_{k'}=Z \otimes k' \subset (X_1 \otimes k') \times (X_2 \otimes k')$
and there exists a local coordinate system at $a$ in
$(X_1 \otimes k') \times (X_2 \otimes k')$ such that the regular
function
$$
t \mapsto \tau_n((Z_{k'})_{t}, (\alpha_{k'})_t)
$$
is defined by the polynomials with the same coefficients as the function
$$
t \mapsto \tau_n(Z_t, \alpha_t).
$$
In model-theoretic terms, this observation implies that once a point
$a \in M^2$, a local coordinate system at $a$, and a family of
branches $\alpha$ of $Z$ are fixed, the slope $\tau_n(Z_t, \alpha_t)$
is definable in the language of fields over $t$.

Further if $X_1 \times \ldots \times X_n$ is a product of
$k$-varieties, we denote the natural projections
$$
p_{i_1 \ldots i_k}: X_1 \times \ldots \times X_n \to X_{i_1} \times
\ldots \times X_{i_k}
$$
denote the natural projections. Although the notion of the composition
of correspondences is standard, we reintroduce it here to fix
conventions.

\begin{defn}[Composition of curves]
  \label{compo}
  Let $Z \subset X_1 \times X_2 \times T$,
  $W \subset X_2 \times X_3 \times S$ be families of curves, and let
  $p_{i_1\ldots i_k}$ denote projections on products of the factors of
  the space $X_1 \times X_2 \times X_3 \times T \times S$. Define the
  \emph{family $W \circ Z$ of compositions of curves from the families
    $W$ and $Z$} to be
  $$
  p_{1345}(p_{124}^{-1}(Z) \cap p_{235}^{-1}(W))
  $$
  in $X_1 \times X_3 \times T \times S$. Clearly, if $Z, W$ are
  definable then so is $Z\circ W$; on the level of points:
  $$
  W \circ Z =\suchthat{ (x,z,t,s) \in M^2 \times T \times S}{ \exists u\
    (x,y,u,t) \in Z \textrm{ and } (y,z,v,s) \in W }.
  $$
  If for all $t \in T, s \in S$ all irreducible components of
  $Z_t, W_s$ project dominantly on $X_1, X_2$, respectively, then
  $W \circ Z$ is a family of curves parametrized by $T\times S$.
  
  We denote by $Z^{-1}$ the image of $Z$ under the morphism
  $X_1 \times X_2 \times T \to X_2 \times X_1 \times T$ that permutes
  the factors $X_1$ and $X_2$, in both geometric and definable
  contexts. We regard the above definitions as applicable to
  individual curves $Z, W$ by putting $T=S$ to be a point.
\end{defn}

\begin{rem}
  If $Z, W$ are families of pure-dimensional curves such that for all
  $t \in T, s \in S$ all irreducible components of $Z_t, W_s$ project
  dominantly on $X_1, X_2$, respectively, then $W \circ Z$ is a family
  of pure-dimensional curves. 
\end{rem}

\begin{prop}
  \label{compo-slope}
  Let $Z \subset X_1 \times X_2 \times T$ and
  $W \subset X_2 \times X_3 \times S$ be families of pure-dimensional
  curves, let $\alpha$, $\beta$ be families of branches of $Z$, $W$ at
  $a=(a_1,a_2) \in Z, b=(b_1,b_2) \in W$, respectively,
  $a_2=b_1$. Then there exists a family of branches
  $\beta \circ \alpha$ of $W \circ Z$ at $(a_1, b_2)$ such that for
  all $t \in T, s \in S$ and for all $n > 0$
  $$
  \tau_n(W_s \circ Z_t, (\beta \circ \alpha)_{(t,s)}) = \tau_n(W_s,
  \beta_s) \circ \tau_n(Z_t, \alpha_t)
  $$
  where the operation ``$\circ$'' in the right hand side expression is 
  composition of truncated polynomials.
\end{prop}

\begin{proof}
  The proof consists essentially in unraveling the definitions. The
  choice of coordinate systems induces the isomorphisms
  $$
  \what{\Oo_{X_1 \times X_2, a}} \cong k[[x,y]] \qquad
  \what{\Oo_{X_2 \times X_3, b}} \cong k[[y,z]] 
  $$
  If the family of branches $\alpha$ is given Zariski locally around
  $t \in T$ by an equation $y - f$, $f \in x\,\Oo_{T,t}[[x]]$, and
  $\beta$ is given by $z - g, g \in y\,\Oo_{S,s}[[y]]$, then let the
  family of branches $\beta \circ \alpha$ be given by
  $z - g \circ f, g \circ f \in (\Oo_{T,t} \otimes
  \Oo_{S,s})[[x]]$. Note that the composition $g \circ f$ of the
  formal power series makes sense and has a zero constant term, since
  both $f$ and $g$ have this property.
  
  Now let $h_Z, h_W$ be generators of the kernels of the maps
  $\Oo_{X_1 \times X_2 \times T,(a,t)} \to \Oo_{Z,(a,t)}$,
  $\Oo_{X_2 \times X_3,(b,s)} \to \Oo_{W,(b,s)}$, respectively, then
  $y-f$ divides $h_X$ and $z - g$ divides $h_Y$. The germ of
  $W \circ Z$ around $(a_1, b_2, t,s)$ by Definition~\ref{compo} is
  cut out by the ideal $(h_X,h_Y) \cap k[[x,z]]$, and in order to show
  that $\beta \circ \alpha$ is a family of branches of $Y\circ X$ at
  this point, we need to check that $(z - g\circ f)$ contains
  $(h_X, h_Y) \cap (\Oo_{T,t} \otimes \Oo_{S,s})[[x,z]]$, and for that
  it would suffice to check that
  $$
  (z - g \circ f) = I := (y - f, z - g) \cap (\Oo_{T,t} \otimes
  \Oo_{S,s})[[x,z]].
  $$
  Indeed, it is straightforward to check that for any $n > 0$
  $$
  (z - g\circ f) = I_n := I/(x^n,z^n)
  $$
  and since $I$ is the inverse limit of $I_n$, it follows that
  $(z - g \circ f) = I$.
\end{proof}

\begin{defn}[Point-wise addition of curves]
  \label{sum}
  Let $G$ be a 1-dimensional algebraic group, and let
  $X \subset G^2 \times T$, $Y \subset G^2 \times S$ be families of
  curves. Let $a: G \times G \to G$ be the group law, let
  $\Gamma_a \subset G^3$ be its graph, and denote $p_{i_1 \ldots i_k}$
  projections of $G \times G \times G \times G \times T \times S$ on
  the products of factors. We define the family of curves
  \emph{$X + Y$ of sums of elements of the families $X$ and $Y$} to be
  $$
  p_{1456}(p_{234}^{-1}(\Gamma_a) \cap p_{124}^{-1}(X)
  \cap p_{135}^{-1}(Y))
  $$
  in $G^2 \times T \times S$. Clearly, if $X, Y$ are definable then so
  is the family $X + Y$; on the level of points:
  $$
  X + Y := \{ (a,b+c, t, s) \mid (a,b,t) \in X, (a,c,s) \in Y \}.
  $$      
\end{defn}

\begin{rem}
  If $X$ and $Y$ are families of pure-dimensional curves then so is
  $X + Y$.  It may seem from the notation above that $G$ is supposed
  to be commutative, even though the definition applies even if this is
  not the case. In this paper we will only consider the operation
  ``+'' for curves inside groups with the commutative connected
  component of the identity.
\end{rem}

Let $G$ be a one-dimensional algebraic group, then the \emph{formal
  group law of $G$} is defined as the image of the topological
generator of $k[[x]]\cong \what{\Oo_{G,e}}$ under the morphism
$\what{\Oo_{G,e}} \to \what{\Oo_{G,e} \otimes \Oo_{G,e}} \cong
k[[x,y]]$ induced by the group operation morhpism. The truncation to
first order of a one-dimensional formal group law is $x + y$
(see, for example, \cite[I.2.4]{jantzen2007representations}).

\begin{prop}
  \label{sum-slope}

  Let $G$ be a one-dimensional algebraic group over an algebraically
  closed field $k$. Let $F$ be the formal group law of $G$, and let
  $F_n$ be its $n$-th order truncation. Let
  $X \subset G \times G \times T$, $Y \subset G \times G \times S$ be
  families of pure-dimensional curves and let $\alpha, \beta$ be
  families of branches at $a=(a_0,a_1), b=(b_0,b_1)$, where $b_0=a_0$,
  respectively. Then there exists a family of branches
  $\alpha + \beta$ of $X+Y$ at $(a_0, a_1+b_1)$ such that
  $$
  \tau_n(X_t + Y_s, \alpha_t + \beta_s)(x) = F_n(\tau_n(X_t,
  \alpha_t)(x) ,\tau_n(Y_s, \beta_s)(x))
  $$
  if $\tau_n(X_t, \alpha_t)$ and $\tau_n(Y_s, \beta_s)$ are
  defined. In particular, if $n=1$,
  $$
  \tau_1(X_t + Y_s, \alpha_t + \beta_s)(x) = \tau_1(X_t,
  \alpha_t) + \tau_1(Y_s, \beta_s).
  $$
\end{prop}

\begin{proof}
  As in the proof of Proposition~\ref{compo-slope}, this statement
  follows from the unfolding of the definitions. Reasoning locally,
  assume that $\alpha$ is cut out by the equation $y - f$, $\beta$ is
  cut out by $z - g$. Then $\alpha + \beta$ is cut out by
  $u - F(f(x), g(x))$. Checking that this is indeed a branch of
  $X + Y$ is straightforward and we leave it to the reader.
\end{proof}

\subsection{Flat families and definability of tangency}
\label{flat-families}

In the present section we prove the main technical result of the
paper, Proposition~\ref{tangency-intersection}, characterizing the
tangency of two generic elements of two families of curves in terms of
properties of the families definable in the reduct $\Mm$. While we do
not give a full definable characterization of the tangency, we prove
a standard weakening of this result, which, as we will see in the
concluding section of the paper, is sufficient for our needs.

The key preliminary step is the observation that if
$X \subset M^2 \times T$, $Y \subset M^2 \times S$ are families of
pure-dimensional curves and $M, T, S$ are smooth, then the `family of
intersections' $X \times_{M^2} Y \to T \times S$ is flat if
restricted to the open subset of $T \times S$ over which it has finite
fibres. We refer the reader to any standard exposition of flatness,
such as \cite[I.\S 2]{milne}, for details, and quickly recall some of
the key facts. All schemes in this section are assumed Noetherian, and
by varieties we mean schemes of finite type over an algebraically
closed field $k$. We identify closed scheme-theoretic points of
varieties with geometric points (that is, morphisms $\Spec k \to X$).

First recall that a morphism $f: X \to Y$ is flat if all local rings
$\Oo_{X,x}$ are flat $\Oo_{Y,f(x)}$-modules. In particular, flatness
can be checked Zariski locally on the source: if
$X = \bigcup\limits_{i=1}^n O_i$ is an open cover and $O_i$ is flat
over $Y$  for all $i$ then $X$ is flat over $Y$ {\cite[\href{https://stacks.math.columbia.edu/tag/01U5}{Tag 01U5}]{stacks}}.

\begin{fact}[{Generic Flatness, \cite[IV.6.10-11]{sga1}}]
  \label{generic-flatness}

  Let $Y$ be an integral scheme and let $f: X \to Y$ be a dominant
  morphism of finite type. Then there exist open subsets
  $O \subset X, U \subset Y$ such that
  $f|_O: O \to Y, f|_{f^{-1}(U)}: f^{-1}(U) \to U$ are flat.
\end{fact}

\begin{fact}[{\cite[Propositions~I.2.4-5]{milne}},
  {\cite[\href{https://stacks.math.columbia.edu/tag/05BC}{Tag
      05BC}]{stacks}}, {\cite[\href{https://stacks.math.columbia.edu/tag/02KB}{Tag
    02KB}]{stacks}} ]
  \label{flat} 

  \noindent\begin{enumerate}[label=(\roman*)]
  \item
    \label{open}
    an open immersion is flat;

  \item
    \label{comp-flat}
    a composition of flat morphisms is flat;

  \item 
    \label{base-change}
    let $X \to Y$ be a flat morphism and let $Z \to Y$ be a
    morphism. Then $X \times_Y Z \to Z$ is flat;

  \item
    \label{criterion}
    let $B$ be a flat $A$-algebra and consider $b \in B$. If the
    image of $b$ in $B/\m B$ is not a zero divisor for any maximal
    ideal $\m$ of $A$ then $B/(b)$ is a flat $A$-algebra;

  \item
    \label{locfree}
    a finite morphism $f: X \to Y$ is flat if and only if $f$ is a
    locally free morphism, that is, if $f_* \Oo_X$ is a locally free
    $\Oo_Y$-module;
    
  \item
    \label{completion}
    if $A$ is an algebra and $I \subset A$ is an ideal, then the
    completion $\varprojlim A/I^n$ is flat over $A$.
  \end{enumerate}
\end{fact}

\begin{lemma}
  \label{pure-dim-flat}
  Assume that the total space $X$ of a family of pure-dimensional
  curves $X \subset M^2 \times T$ is a variety. Then $X$ flat over
  $T$.
\end{lemma}

\begin{proof}
  By definition of a family of pure-dimensional curves $X$ is open in
  $\bar X$, so by Fact~\ref{flat}\ref{open} and \ref{comp-flat}
  suffices to show flatness of $\bar X$ over $T$. Passing to a cover
  of $M^2 \times T$ by affine opens $O_i = \Spec B_i$ suffices to show
  flatness of $\bar X \cap O_i$ over $T$ for all $i$. But this
  immediately follows from the definition of a family of
  pure-dimensional curves and Fact~\ref{flat}\ref{criterion}.
\end{proof}

\begin{lemma}
  \label{int-flat} 
  Let $M$ be a smooth algebraic curve, $T,S$ be smooth varieties, 
  $X \subset M^2 \times T$, $Y \subset M^2 \times S$ be families of
  pure-dimensional curves, and assume that $X, Y$ are varieties. Let
  $U \subset T \times S$ be the set of points $u$ such that
  $\dim (X \times_{M^2} Y)_u = 0$, let
  $Z = X \times_{M^2} Y \cap p^{-1}(U)$, where $p$ is the projection
  onto $T \times S$. Then the restriction $p: Z \to U$ is flat.
\end{lemma}

\begin{proof}
  By Lemma~\ref{pure-dim-flat}, $X$ is flat over $T$.  Since $M, T, S$
  are smooth and since regular local rings are unique factorization
  domains, and $Y$ is pure-dimensional, $Y$ is cut out in
  $M^2 \times S$ by a principal ideal sheaf (see, for example,
  \cite[\S 19, Theorems~47, 48]{matsu-ca}). Since $X \to T$ is flat,
  by Fact~\ref{flat}\ref{base-change}
  $X \times S \cong X \times_T (T \times S) \to T \times S$ is flat
  too. Since $X$ is pure-dimensional, the natural closed embedding
  $X \times_{M^2} Y \to X \times S$ is also cut out by a principal
  ideal sheaf $\Ii$, passing by virtue of Fact~\ref{flat}\ref{open} to
  a cover of $X \times S$ by affine opens $O_i=\Spec B_i$ and applying
  Fact~\ref{flat}\ref{criterion} to the algebras $B_i$, this closed
  subscheme is flat precisely over the complement of the subvariety of
  $T \times S$ consisting of the points $u$ such that the local
  generator of $\Ii$ does not vanish on an irreducible component of
  the fibre $(X \times S)_u$. In other words, it is flat over the open
  subset of points $u \in T \times S$ such that $(X \times_{M^2} Y)_u$
  is zero-dimensional.
\end{proof}

Recall that a morphism $f: X \to Y$ is called \emph{finite} if for any
affine open $U=\Spec R \subset Y$ the inverse $f^{-1}(U)=\Spec S$ is
affine and $S$ is a finite $R$-module. A morphism is finite if and
only if it is quasi-finite and proper (\cite[EGA~IV.18.12.4]{ega}).

\begin{lemma}
  \label{sc} 
  Let $f: X \to Y$ be a quasi-finite morphism of schemes over an
  algebraically closed $k$. Consider the functions
  $$
  \begin{array}{ll}
  m: Y(k) \to \Z & y \mapsto \#(f^{-1}(y)), \\
  w: X(k) \to \Z & x \mapsto \dim_k \Oo_{X,x} \otimes k(f(x)).\\
  \end{array}
  $$  
  Then
  \begin{enumerate}[label=(\roman*)]
  \item
    \label{i}
    $w(x) = \dim_{k(x)} \what{\Oo_{X,x}} \otimes k(f(x))$ where
    $ \what{\Oo_{X,x}}$ is the completion
    $ \varprojlim \Oo_{X,x}/I^n$ for any ideal $I \subset \Oo_{X,x}$;
  \item
    \label{ii}
    assume that $f$ is flat. Then $m$ is lower semi-continuous and $w$
    is upper semi-continuous, that is, the lower level sets of $m$ and
    the upper level sets of $w$
    $$
    \suchthat{y \in Y}{m(y) \leq n}
    \textrm{ and } \suchthat{x \in X}{w(x) \geq n}
    $$
    are closed; $y \mapsto \sum_{x \in f^{-1}(y)} w(x)$ is lower
    semi-continuous  and is locally constant if $f$ is finite.
  \end{enumerate}  
\end{lemma}

\begin{proof}
  Let $J$, resp. $\hat J$, be the ideal of $\Oo_{X,x}$,
  resp. $\what{\Oo_{X,x}}$, generated by the image of the maximal ideal of
  $\Oo_{Y,f(x)}$, then
  $\Oo_{X,x} \otimes k(f(x)) \cong \Oo_{X,x} / J$. We have a sequence
  of isomorphisms
  $$
  \Oo_{X,x} \otimes k(f(x)) \cong \Oo_{X,x} / J \cong
  \what{\Oo_{X,x}} / \hat J \cong \what{\Oo_{X,x}} \otimes k(f(x))
  $$
  where the first and the third one are tautological, and the second
  morphism is an isomorphism because $\hat J$ is the completion of $J$
  in the $I$-adic topology. This proves claim~\ref{i}.

  That $m$ is lower semi-continuous follows from
  \cite[EGA~IV.3.15.5.1(i)]{ega} and the fact that flat morphisms of
  finite type are universally open \cite[EGA~IV.2.4.6]{ega}.  Upper
  semi-continuity of $w$ follows from
  \cite[\href{https://stacks.math.columbia.edu/tag/0F3D}{Tag
    0F3D3}]{stacks} and
  \cite[\href{https://stacks.math.columbia.edu/tag/0F3I}{Tag
    0F3I}]{stacks} (note that the definition of $w$ from the statement
  of the Lemma and one from \cite{stacks} coincide on the closed
  scheme theoretic points of a variety over an algebraically closed
  field).  Lower semi-continuity of $\sum_{x \in f^-1(u)} w(u)$
  follows from
  \cite[\href{https://stacks.math.columbia.edu/tag/0F3J}{Tag
    0F3J}]{stacks}, that it is locally constant if $f$ is finite
  follows from the definition of a weighting
  \cite[\href{https://stacks.math.columbia.edu/tag/0F3A}{Tag
    0F3A}]{stacks}.

\end{proof}

We can now formulate our main technical result. Roughly, it states that,
in suitably chosen families of curves tangency of two curves is
witnessed by a lower number of intersection points:

\begin{prop}
  \label{tangency-intersection} 
  Keep notation and assumptions of Lemma~\ref{int-flat} and assume
  further that there exists $a \in M^2$ such that $X_t, Y_s$ pass
  through $a$ for all $t \in T, s \in S$. Let $\alpha, \beta$ be
  families of branches at $a$ of $X$, $Y$ respectively, such that for
  all $t \in T, s \in S$ the slopes of $\alpha_t, \beta_s$ are
  defined. Define
  $$
  \begin{array}{rcl}
    n_{\max} & = & \max \{ n \mid
              \forall (t,s) \in U(k)\ 
              \tau_n(X_t, \alpha_t) =            
              \tau_n(Y_s, \beta_s)  \}\\
    m_{\max} & = & \max\limits_{(t,s) \in U(k)} \#(X_t \cap  Y_s) \\
  \end{array}
  $$
  Then
  \begin{multline*}
  \suchthat{ (t,s) \in U(k) }{\tau_{n_{\max}+1}(X_t, \alpha_t) =
    \tau_{n_{\max}+1}(Y_s, \beta_s) }  \subseteq \\
  \subseteq
  \suchthat{(t,s) \in U(k)}{\#(X_t\cap Y_s) < m_{\max}} 
  \end{multline*}
\end{prop}

\begin{proof}
  Consider $Z = X \times_{M^2} Y \cap M^2 \times U$, let
  $q: Z \to M^2$ be the natural projection and let
  $Z = \bigcup_{i=0}^n Z_i$ be the decomposition into irreducible
  components where $Z_0=q^{-1}(a)$. We will first show that whenever
  $p^{-1}(u) \cap Z_0 \cap Z_i \neq \emptyset$ for some $i \neq 0$ we
  have $\#p^{-1}(u) < \max_{u \in U}\#p^{-1}(u)$. In order to do that
  we will show that the function
  $u \mapsto \#p^{-1}(u) \cap (\cup_{i\neq 0} Z_i)$ is lower
  semi-continuous.

  The projection $p: Z \to U$ is flat by Lemma~\ref{int-flat}. By
  \cite[\href{https://stacks.math.columbia.edu/tag/04PW}{Tag
    04PW}]{stacks} the closed embedding $Z_{\red} \to Z$, where
  $Z_{\red}$ is $Z$ endowed with the canonical reduced structure, is
  flat.  Since the invariant we are interested in does not depend on
  the scheme structure, by Fact~\ref{flat}\ref{comp-flat} we may
  assume $Z$ reduced. Furthermore, there exists an open embedding
  $Z \hookrightarrow \bar Z$ where $\bar Z$ is flat and finite over
  $U$. Indeed, let $\bar M$ be a smooth proper algebraic curve that
  contains $M$ as a dense subset and let $\bar X, \bar Y$ be the
  closures of $X, Y$ in $\bar M^2 \times T$ and $\bar M^2 \times
  S$. Let
  $\bar Z = \bar X \times_{\bar M^2} \bar Y \cap \bar M^2 \times U$,
  let $\bar p$ be the natural projection on $U$ and denote $\hat p$
  its restriction to $\hat Z=\cup_{i\neq 0} \bar Z_i$, where
  $\bar Z_i$ is the irreducible component of $\bar Z$ that contains
  $Z_i$ for each $i$. By Lemma~\ref{int-flat} $\bar p$ is flat.

  By Fact~\ref{flat}\ref{locfree} the morphism $\bar p$ is locally
  free.  It is readily seen that $(\bar p)_* \Oo_{\hat Z}$ is locally
  free of rank one less than the rank of $(\bar p)_* \Oo_{\bar
    Z}$. Indeed, if $W = \Spec A \subset U$ is an affine open such
  that $(\bar p)^{-1}(W) = \Spec B$ and $B$ is a free $A$-module, we
  have that $B \cong B/{\mathfrak p} \oplus B/{\mathfrak q}$ as
  $A$-module, where $\mathfrak p, \mathfrak q \subset B$ are ideals
  cutting out
  $Z_0 \cap (\bar p)^{-1}(W), \hat Z \cap (\bar p)^{-1}(W)$,
  respectively. Since $Z_0 \cong U$, in particular
  $Z_0 \cap (\bar p)^{-1}(W) \cong W$, and we have that
  $B/{\mathfrak p} \cong A$, so $B/{\mathfrak q}$ is free. By
  Fact~\ref{flat}\ref{locfree} again $\hat p$ is flat, and by
  Fact~\ref{flat}\ref{open} its restriction to
  $\hat Z \cap Z = \cup_{i \neq 0} Z_i$ is flat.  We deduce by
  Lemma~\ref{sc}\ref{ii} that the function
  $u \mapsto \#p^{-1}(u) \cap (\cup_{i \neq 0} Z_i)$ is lower
  semi-continuous.

  Note that while $Z_0$ may have non-trivial scheme-theoretic
  structure, the restriction $p|_{Z_0}: (Z_0)_{\red} \to U$ is a
  homeomorphism, so denote $r: U \to Z_0$ its set-theoretic
  inverse. Let $w: Z \to \Z, w(z) = \dim_k \Oo_{Z,z} \otimes
  k(p(z))$. We claim that $w$ is constant on the open set
  $Z'=Z_0 \wo \bigcup_{i\neq 0} Z_i$. By Fact~\ref{flat}\ref{open}
  $Z'$ is flat over $U$ and by Fact~\ref{flat}\ref{base-change} it is
  flat over the open $p(Z') \subset U$. The restriction
  $p|_{Z'}: Z' \to p(Z')$ is still a homeomorphism; we will show that
  it is a finite morphism.

  Note that since $U$ is dense in $T \times S$ it is integral, and
  since $Z'$ is dense in $U$, it is integral too. The scheme $Z'$ is
  of finite type over a field, so clearly quasi-compact, and $p$ is
  clearly separated (for example, because it is affine), so Zariski's
  Main Theorem (see \cite[Ch.~I, Theorem~1.8]{milne}) can be applied
  to $p$. Therefore, $p$ factors into a composition of an open
  embedding $i: Z' \to Z''$ and a finite morphism $p': Z'' \to
  p(Z')$. Since $\tilde p=p_{\red}|_{Z'}: Z'_{\red} \to p(Z')$ is an
  isomorphism, the morphism
  $p''=i_{\red} \circ \tilde p^{-1}\circ p'_{\red} : Z''_{\red} \to
  Z''_{\red}$ restricts to the identity morphism on $Z'_{\red}$. If
  $p|_{Z'}: Z' \to p(Z')$ is not finite then the open embedding $i$ is
  an isomorphism and $p''$ is not an isomorphism.  Since passing to a
  closed subscheme preserves finiteness, we may assume $Z'$ to be
  dense in $Z''$. The subset of $Z''_{\red}$ where $p''$ and the
  identity morphism coincide is closed and contains $Z'$, so must be
  the whole of $Z''$, which in turn contradicts $p''$ not being an
  isomorphism.

  Now by Lemma~\ref{sc}\ref{ii} $w$ is upper semi-continuous on $Z(k)$
  and in particular on $Z_0$, but since $Z'$ is flat and finite over
  $U$, $w$ is constant on $Z'(k)$. Therefore, $w$ takes the value
  $w_{\min,0} = \min_{x \in Z_0} w(x)$ on the latter, and if
  $w(r(u)) > w_{\min,0}$ for some $u \in U$ then
  $r(u) \in Z_i \cap Z_0$ for some $i \neq 0$ and therefore
  $\#p^{-1}(u)$ is not maximal. It follows that
  $$
  \suchthat{ (t,s) \in U(k)}{ w(r(u)) > w_{\min,0}}  \subseteq 
  \suchthat{ (t,s) \in U(k)}{\#(X_t \cap Y_s) < m_{\max}}.
  $$
  It is left to prove that
  $$
  w(r(t,s)) > w_{\min,0} \textrm{ for all }t,s \textrm{ such that }
  \tau_{n_{\max}+1}(X_t, \alpha_t) = \tau_{n_{\max}+1}(Y_s, \beta_s).
  $$
  To establish this, it is enough to prove the statement on an affine
  Zariski open subset $\Spec R \subset U \times M^2$ intersecting
  $Z_0$ non-trivially. Let $f,g \in R$ be the equations of
  $X \times S \cap \Spec R, Y \times T \cap \Spec R$,
  respectively. Let $I \subset R$ be the ideal of functions that
  vanish on $q^{-1}(a)$ and let $\what R = \varprojlim R/I^n$.
  
  Let $f=f_1 \cdot \ldots \cdot f_N$, $g = g_1\cdot \ldots \cdot g_K$
  be decompositions into pairwise coprime factors in $\what R$ and let
  $f_\alpha$ and $g_\beta$ be those factors that are local equations
  of $\alpha, \beta$. Apply the Chinese Remainder Theorem (see
  \cite[Ch.~9, ex.~9]{atiyahmac-ca}) twice: first, to $\what{R}/(f,g)$
  to get the decomposition
  $$
  \what{R}/(f,g) = \bigoplus_{i=1}^N  \what{R}/(f_i, g),
  $$
  second, to each direct summand $\what{R}/(f_i,g)$ to get
  $$
  \what{R}/(f,g) = \bigoplus_{i=1}^N \bigoplus_{j=1}^K  \what{R}/(f_i,
  g_j). 
  $$
  Both applications are justified: since $f_i$ are pairwise coprime in
  $\what{R}$ (that is, $(f_i) + (f_k) = \what{R}$ for $i\neq k$), the
  ideals $(f_i,g)$ are pairwise coprime in $\what{R}/(g)$,
  similarly, for each $i$, the ideals $(f_i, g_j)$ are pairwise
  coprime in $\what{R}/(f_i,g)$.

  Tensoring with $k(u)$ and applying Lemma~\ref{sc}\ref{i} we
  get
  $$  
  w(r(u)) = \sum_{i=1}^N \sum_{j=1}^K \dim_{k(u)} \what{R}/(f_i,g_j)\otimes k(u).
  $$
  Therefore, if $u=(t,s) \in U(k)$ and
  $\tau_{n_{\max}+1}(X_t, \alpha_t) = \tau_{n_{\max}+1}(Y_s, \beta_s)$
  then by the remark after Definition~\ref{relslope-definition}
  $$
  \dim_k \what{R}/(f_\alpha, g_\beta) \otimes k(u)
  $$
  takes a value strictly greater than the minimum it achieves on
  $U(k)$. By Fact~\ref{flat}\ref{completion} $\what R$ is a flat
  $R$-algebra, and by applying Fact~\ref{flat}\ref{criterion} twice,
  as in the proof of Lemma~\ref{int-flat}, $\Spec \what R/(f_i,g_j)$
  is flat over $U$ for all $i,j$. Since by Lemma~\ref{sc}\ref{ii} for
  each pair of prime factors $f_i,g_j$ the value
  $$
  \dim_{k} \what R/(f_i, g_j) \otimes k(u)
  $$
  is upper semicontinuous in $u$, it follows that
  $w(r(t,s)) > w_{\min,0}$ as soon as slopes of order $(n_{\max} + 1)$
  of $\alpha_t$ and $\beta_s$ coincide.
\end{proof}

\section{Interpretation of the field}
\label{interp-field}

In the preset section we tie together the results obtained above to produce the main result of the paper. We start with some additional technicalities and reductions. 

\subsection{The group configuration}\label{S:gc}

In stable theories -- a model theoretic framework encompassing all
structures considered in the present work -- certain combinatorial
configurations of elements are known to exist only in the presence of
an interpretable group or -- in a more restrictive setting -- an
interpretable field. It is by constructing such configurations --
using the ``definable intersection theory'' developed in the previous
sections -- that the main theorem of the present paper is proved.

Before describing these configurations in more detail we need some
model theoretic preliminaries. As in Section~\ref{background}, we will
specialise the definitions to the setting in which they will be
used. As above, we will be working in the full Zariski structure $\M$
on an algebraic curve $M$ over an algebraically closed field $k$. We
will be mostly concerned with a structure $\Mm:=(M,X)$ where
$X\subseteq M^2\times T\subseteq M^{2+l}$ is the total space of an
ample family. Throughout the text by \emph{definable} we mean
`definable with parameters'.

\begin{defn}
	\begin{enumerate}
		\item 	If $D$ is an $\Mm$-definable set we let $\dim(D):=\dim(\mathrm{cl}(D))$ where $\mathrm{cl}$ denotes the Zariski closure of $D$. 
		\item If $A$ is a set of parameters and $a\in M^l$ we denote $\dim(a/A):=\min\{\dim(D): a\in D\}$ where $D$ ranges over all subsets of $D^l$ $\Mm$-definable over $A$. 
		\item We say that $a\in M$ is $\Mm$-algebraic over $A$ if $\dim(a/A)=0$. We denote $\acl_\Mm(A):=\{a\in M: \dim(a/A)=0 \}$. 
		\item We say that $a$ is $\Mm$-generic in $D$ over $A$ if $\dim(a/A)=\dim(D)$. 
		\item We say that $a$ is $\Mm$-independent from $B$ over $A$ if $\dim(a/A)=\dim(a/AB)$. 
	\end{enumerate}
\end{defn}

\begin{rem}
  Note that $\dim(D)$ is (by definition) the same as the
  algebro-geometric dimension of $D$. This implies that
  $\dim_\Mm(a/A)\ge \dim_\M(a/A)$, but there is no need for equality
  to hold. In particular $\acl_\Mm(A)\subseteq \acl(A)$ where on the
  right hand side $\acl$ is the field-theoretic algebraic closure in
  $k$. It follows that $\M$-independence (which coincides with the
  field-theoretic notion of algebraic independence) implies
  $\Mm$-independence, but not necessarily the other way around.
\end{rem}

\begin{defn}
  An infinite set $D$ definable (or interpretable) in $\Mm$ is
  \emph{strongly minimal } if every $\Mm$-definable subset of $D$ is
  finite or cofinite.
\end{defn}

It is an easy exercise to verify that if $D$ is an $\Mm$-definable set and $\dim(D)>1$ then there is a projection $\pi:D\to M^{\dim(D)-1}$ and an open $U\subseteq M^{\dim(D)-1}$ such that $\pi^{-1}(u)\cap D$ is infinite for all $u\in U$. In particular, $D$ is strongly minimal only if $\dim(D)=1$. Thus, $D$ is strongly minimal if and only if it is one-dimensional and cannot be written as the disjoint union of two one-dimensional $\Mm$-definable subsets. We say that $\Mm$ is strongly minimal if $M$ is (as an $\Mm$-definable set). 

\begin{rem}
  An $\Mm$-definable set $D$ may be strongly minimal with
  respect to the structure $\Mm$ but not with respect to the
  strucutre $\M$. 
\end{rem}

As we will see below, we can easily reduce the proof of our main result to the case where $\Mm$ is strongly minimal. Under this additional assumption we can finally introduce the group configuration: 

\begin{defn}[Group configuration]\label{def-gconf}
	\label{def:grconf}
	
	Let $\Mm$ be as above, and assume that it is strongly minimal.
	The set $\set{a,b,c,x,y,z}$ of tuples
		$$
	\groupconf{a}{b}{c}{x}{y}{z}
	$$
	 is a \emph{group
		configuration} if there exists an integer $n$ such that
	\begin{itemize}
        \item all elements of the diagram are pairwise independent and
          $\dim(a,b,c,x,y,z)=2n+1$;
        \item $\dim a = \dim b = \dim c=n$,
          $\dim x = \dim y = \dim z = 1$;
        \item all triples of tuples lying on the same line are
          dependent, and moreover,\linebreak $\dim(a,b,c)=2n$, 
          $\dim(a,x,y)=\dim(b,z,y)=\dim(c,x,z)=n+1$;
	\end{itemize}
	Two group configurations $G_1, G_2$ are called \emph{inter-algebraic}
	if for any pair of tuples $a \in G_1, a' \in G_2$ in the
	corresponding vertices, $\acl_\Mm(a) = \acl_\Mm(a')$.
\end{defn}

Assume that $G$ is a connected $\Mm$-definable group acting
transitively on a strongly minimal definable set $X$, then one can
construct a group configuration as follows: let $g,h$ be independent
generics in $G$ and let $u$ be a generic of $X$ (we will justify the
assumption that such generics exist later on), then
$\tuple{g,h,g\cdot h, u, g \cdot u, g \cdot h \cdot u}$ is a group
configuration (associated with the action of $G$ on $X$). Below (Lemma
\ref{many-slopes}) we show that, for a suitably constructed
$\Mm$-definable family of curves passing through a fixed point, the
set of $n$-th slopes of curves in the family coincide for some $n$
with a one-dimensional algebraic group, $H$ (viewed as acting on
itself by multiplication). Proposition~\ref{tangency-intersection}
will then allow us to `pull back' a group configuration (in $\M$)
associated with this group $H$ into a group configuration in
$\Mm$. This will, essentially, finish the proof since:

\begin{fact}[Hrushovski]
  \label{grconf} 
	
  Let $M$ be a strongly minimal structure and let
  \linebreak$G_1 = \tuple{a,b,c,x,y,z}$ be a group configuration. Then
  there exists a definable group $G$ acting transitively on a strongly
  minimal set $X$. 
\end{fact}

This follows from the main theorem of \cite{bouscaren-group} and the
fact that infinitely definable groups in stable theories are
intersections of definable groups (see, for example,
\cite[Theorem~5.18]{poizat}) and the fact that any group definable in
an algebraically closed field is (definably isomorphic to) an
algebraic group (see \cite[Theorem~4.13]{poizat}). The original proofs
of these statements
are contained in \cite{hru-phd}.\\

To construct a field we will have to work a little harder. First, 

\begin{defn}
		A group configuration $\tuple{a_1,a_2,a_3,x,y,z}$ is \emph{minimal} if 
		\begin{align*}
			&	\acl_\Mm(\mathrm{Cb}(x,y)/a_1)=\acl_\Mm(a_1), \\
		 	& 	\acl_\Mm(\mathrm{Cb}(y,z)/a_2)=\acl_\Mm(a_2), \\
		 	& 	\acl_\Mm(\mathrm{Cb}(x,z)/a_3)=\acl_\Mm(a_3).
		\end{align*}
\end{defn}

\begin{rem} 
	We will not go into the definition of canonical bases (see, e.g., \cite[p.19]{pillay1996geometric}), but for the benefit of readers unfamiliar with this model theoretic notion we mention that: 
	\begin{enumerate}
		\item The minimality condition is readily checked to be equivalent to the condition that whenever there are $a'_i\in \acl_\Mm(a_i)$ such that $\tuple{a'_1,a'_2,a'_3,x,y,z}$ is still a group configuration then $a_i\in \acl_\Mm(a_i')$ for all $i=1,2,3$.
		\item By dimension considerations it follows from the previous remark that any group configuration $(a_1,a_2,a_3,x,y,z)$ gives rise to a minimal group configuration $(a_1',a_2',a_3',x,y,z)$ with $a_i'\in \acl(a_i)$ for all $i$. In particular, if $\dim(a_i)=1$ for all $i$ then $(a_1,a_2,a_3,x,y,z)$ is a minimal configuration. 
		\item Roughly, $\mathrm{Cb}((x,y)/a)$ is the model theoretic analogue of the field of definition of the locus of $(x,y)$ over $a$. 
		\item For our purposes it will suffice to know that if $X\to T$ is a nearly faithful family of curves (see
        below) then $t$ is (up to inter-algebraicity) a
        canonical base for $x/t$ for any generic point of
        $x$, and if through $x_1,\dots, x_k$ there is only
        one curve $X_t$ in $X$ then $t$ is (up to
        inter-algebraicity) a canonical base of
        $(x_1,\dots, x_k)$.
	\end{enumerate}
\end{rem}


For minimal group configurations we have: 

\begin{fact}[{\cite[Theorem~V.4.5]{pillay1996geometric}}]
	\label{grconf-faithfulness} 	
	If the group configuration in the statement of Fact~\ref{grconf} is, additionally, assumed to be minimal then 
	the action of the group $G$ on $X$ as provided above can be taken to be faithful and this group action has an associated group configuration
	$G_2 = \tuple{g,h,g\cdot h, u, g \cdot u, g \cdot h \cdot u}$ 
	inter-algebraic with $G_1$. In particular, $\dim G=\dim a$.
\end{fact}

This, finally, allows to obtain a field as follows: 

\begin{fact}[Hrushovski, \cite{hru-phd}]
	\label{fieldconf} 
	
	Let $G$ be an $\Mm$-definable group acting transitively and
	faithfully on a strongly minimal set $X$. Then either $\dim(G)=1$ or there exists a
	definable field structure on $X$ and either $\dim(G)=2$ and
	$G \cong \G_a \rtimes \G_m$, or $\dim(G)=3$ and $G \cong \mathrm{PSL}_2$.
\end{fact}

An exposition of the above fact can be found in \cite{poizat}
(Theorem~3.27). Establishing that $G$ is isomorphic to
$\G_a \rtimes \G_m$ or to $\mathrm{PSL}_2$ is the crucial point in the
proof of Fact~\ref{fieldconf}. In the present context, where $G$ and
$X$ are definable in an algebraically closed field (rather, the full
Zariski structure on an algebraic curve) this statement can be
established using a simpler direct algebraic proof.

\subsection{Some standard reductions}

We make some standard simple reductions that will allow us to more
easily use the results obtained in the previous sections as well as
the group and field configurations described above.

\begin{lemma}
	We may assume that $k$ is of infinite transcendence degree (over the prime field). 
\end{lemma}
\begin{proof}
	Let $K\ge k$ be an algebraically closed field extension of infinite transcendence degree. We let $M':=M(K)$, and for any $D$ $\Mm$-definable without parameters we let $D':=D(K)$. We obtain a structure $\Mm':=(M',X')$. By Hilbert's Nullstellensatz and  Chevalley's Theorem (see, e.g., \cite[Corollary 3.2.8]{marker}) $X'$ is ample (if $U\subseteq M$ is an open set witnessing the fact that $X$ is ample then $U'$ witnesses that $X'$ is). Note also that any set $S$ $\M'$-definable without parameters is of the form $D'$ for some $\Mm$-definable set, $D$. 
	
	Assume that a field is interpretable in $\Mm'$. This means that there are $D, E$ $\Mm'$-definable (without parameters) and parameters $\bar a\in K^l$ and $\bar b\in K^n$ such that $E_{\bar b}$ is an equivalence relation (of the correct arity) and such that $D_{\bar a}/E_{\bar b}$ is an infinite field. Let $L_{\bar c}$ and $A_{\bar d}$ be the graphs of multiplication and addition respectively, for $L,A$ $\Mm'$-definable without parameters. 
	
	Consider the set $S$ of all parameters $(\bar x, \bar y, \bar
        z, \bar w)$ such that $E_{\bar y}$ is an equivalence relation
        on $D_{\bar x}$ and $L_{\bar z}, A_{\bar w}$ turn $D_{\bar
          x}/E_{\bar y}$ into an infinite field. We claim that $S$ is
        $\Mm'$-definable without parameters. This is easy since, if
        $C\subseteq M^{r+s}$ is any constructible set then the set
        $\{v\in M^s : |C_v| < \infty \}$ is uniformly bounded, say by $N$. So $C_v$ is infinite if and only if $|C_v|>N$, which is a definable property (of $v$). By Hilbert's Nullstellensatz and  Chevalley's Theorem again $S$ has a point in $k$, meaning that an infinite field is interpretable already in $k$. 
\end{proof}

In model-theoretic terms the above lemma only means that
interpretability of an (infinite) field is a first order property, and
therefore preserved under the passage to elementary substructures. The
most useful -- for our purposes -- property of fields of infinite
transcendence degree is the following consequence of Chevalley's
theorem and the compactness theorem of first-order logic:

\begin{fact}\label{generics}
	If $k$ is of infinite transcendence degree then any $\Mm$-definable set $D$ has generic points over any finite set of parameters $A$. 
\end{fact}


We need the following (weak) version of \cite[Theorem B.1.43]{zg}: 
\begin{fact}
  If $\Mm$ is strongly minimal and not locally modular then there
  exists an ample definable family of curves
  $X\subseteq M^2\times M^l$ with the property that for any $t\in M^l$
  the set $E_t:=\{s\in M^l: |X_t\cap X_s|=\infty\}$ is finite.
\end{fact}

For the purposes of the present paper we call an ample definably
family of curves as above a
\emph{nearly faithful family of curves}. Combined with the (easy) fact that local modularity is preserved under naming parameters (\cite[Remark IV.1.8]{pillay1996geometric}) this gives: 
\begin{fact}\label{F:nearlyfaithful}
  If $\Mm$ is strongly minimal and not locally modular then there
  exists a nearly faithful ample family whose generic members are all
  strongly minimal subsets of $M^2$.
\end{fact}

We can now show: 

\begin{lemma}
	We may assume that $M$ is a smooth curve and that $\Mm$ is strongly minimal. 
\end{lemma}
The proof is well known (see, e.g., \cite[Lemma IV.1.7]{pillay1996geometric} for a much more general statement). For the sake of completeness we outline a simple proof in the present context.  

\begin{proof}
  Clearly, if $S$ is an $\Mm$-definable set, $\mathcal S$ the
  structure with universe $S$ and definable sets $D\cap S^l$ (as $D$
  ranges over $\mathrm{Def}(\Mm)$ and $l$ ranges over
  $\mathbb N_{>0}$), then, if a field $F$ is interpretable in
  $\mathcal S$, then $F$ is already interpretable in $\Mm$. So it will
  suffice to show that there exists a strongly minimal set
  $S\subseteq \Mm$ that is not locally modular with respect to the induced
  structure $\mathcal S$.
	
  Let $X\subseteq M^2\times T$ be an ample family witnessed by an open
  $U\subseteq M^2$. Reducing $U$, if needed, we may assume that if
  $U_i:=p_i(U)$ are the projections of $U$ onto the two $M$ factors,
  then for all $a,a'\in U_1$ the fibre $\pi_1^{-1}(a)\cap U$ is
  infinite and $(\pi_1^{-1}(a)\cap U)=(\pi_1^{-1}(a')\cap U)$ (up to a
  finite set) for all $a,a'\in U_1$. Thus, setting $Y=X\circ X^{-1}$
  we immediately see that $Y$ is ample witnessed by $U_1\times U_1$.
	
  Since $U_1$ is $\Mm$-definable and one-dimensional there exists some
  strongly minimal $M_0\subseteq U_1$. It is clear that
  $Y\cap (M_0^2\times T^2)$ is an ample family in $\Mm$. It would be
  an ample family in the induced structure on $M_0$ if we could
  replace $T$ with some $T_0\subseteq M_0^r$ (for some $r$).
	
  Recall that an irreducible algebraic curve of degree $d$ is uniquely
  determined by any $d+1$ generic enough points on the curve. A
  similar argument would show that any curve in
  $Y\cap (M_0^2\times T^2)$ is uniquely determined (up to a finite
  set) by finitely many generic enough points on the curve. More
  specifically, by the previous fact we may assume that $X$ is nearly
  faithful. So if $X_t$ is a curve and $p\in X_t$ is generic then
  $\dim(p/t)=1$, so $\dim(p/t)<\dim(p/\0)$ and by symmetry,
  $\dim(t/p)<\dim(t/\0)$. Proceeding in a similar way, we get that
  $\dim(t/p_1,\dots, p_r)=0$ for $r=\dim(T)$ and
  $p_1,\dots, p_r\in X_t$ generic enough. So there are only finitely
  many $t'$ such that $p_1,\dots, p_r\in X_{t'}$ and by adding enough
  points we can assure that $X_t$ is uniquely determined (up to a
  finite set) by some $p_1,\dots, p_l\in X_t$. We leave it to the
  reader to verify that the number $l$ of points determining the curve
  can be taken to be independent of $t$.  Since any $X_t$ has
  infinitely many points in $M_0^2$ we get that $M_0$, with the
  induced structure, is, indeed, strongly minimal and not locally
  modular. Replacing $M_0$ with $M_0\cap M_{reg}$, the regular locus
  of $M_0$, the conclusion follows.

\end{proof}

%
%
%

To sum up: 
\begin{cor}\label{reduction}
  To prove Conjecture~\ref{RTC} it suffice to prove: Let $\Mm$ be a
  strongly minimal reduct of the full Zariski structure $\M$ on a
  smooth algebraic curve $M$ over an algebraically closed field $K$ of
  infinite transcendence degree. Then either $\Mm$ is locally modular
  or $\Mm$ interprets a field $K$-definably isomorphic to
  $K$. Moreover, we may assume that the lack of local modularity of
  $\Mm$ is witnessed by a nearly faithful family of curves whose
  generic members are strongly minimal.
\end{cor}

\subsection{Generically unramified projections}
\label{sect:unram}

In order to apply the machinery of slopes and tangency discussed in
Section~\ref{slopes} we need to produce, definably in $\Mm$, large
enough families of curves where these notions are defined and carry
information. Lemma~\ref{projections} below guarantees the former
requirement, namely that for any curve $X \subset M^2$ the slope is
defined on a dense open subset of either $X$ or $X^{-1}$ (uniformly in
parameters). The second requirement is more delicate, as pointed out
for example in the concluding remarks of \cite{marker-pillay}. In more
technical terms, the problem pointed out by Marker and Pillay is that
if the projection $p_2: Z \to M$ is everywhere ramified for a curve
$Z \subset M^2$ (e.g. the curve cut out by the equation $y = x^p$ in
$\A^1 \times \A^1$) then even if $p_2$ is dominant,
$\tau_1(Z, \alpha)=0$ for any branch $\alpha$ at any point of $Z$.
In Lemma~\ref{fr-cancel} and Lemma~\ref{cancel-insep} we develop the
tools allowing us to construct, definably in $\Mm$, curves in $M^2$
whose projections on both factors $M$ are generically unramified.

The following lemma ensures that at least one of the projections on a
factor $M$ of a family of curves is generically \'etale for a general
element of the family, which by Lemma~\ref{etale-branch} implies
existence of slopes for a generic element of the family. The fact that
the support of the module of K\"ahler differentials is closed and
Fact~\ref{generic-flatness} imply that being \'etale and being
unramfied is open on the source. In particular, in order to check
whether a dominant morphism $f: X \to Y$ is \'etale on a dense open
subset of $X$ it suffices to check if $\Omega_{k(X)/k(Y)}=0$, or
equivalently (see \cite[Exercise~6.2.9]{liu}, also
\cite[Lemma~6.1.13]{liu}), if $k(X) \supset k(Y)$ is a separable
extension. We refer the reader to any standard algebraic geometry
reference (e.g.  \cite[Section~6]{liu},
\cite[Section~II.8,IV.2]{hartshorne}) for the details on K\"ahler
differentials and ramification.

\begin{lemma}
  \label{projections} 
  Let $M$ be an irreducible algebraic curve over a field of any
  characteristic. Let $X \subset M^2 \times T$ be a family
  of pure-dimensional curves, and assume that $X$ and $T$ are
  irreducible. Then there exists a dense open subset $T' \subseteq T$
  such that either $p_1: X_t \to M$ or $p_2: X_t \to M$ is generically
  \'etale for all $t \in T'$.
\end{lemma}

\begin{proof}
  Let $\xi$ be the generic point of $T$ in the scheme-theoretic sense.
  Denote $M_\xi = M \otimes k(\xi)$, $X_\xi=X \otimes k(\xi)$. By
  slightly abusing notation, denote $p_1, p_2: X_\xi \to M_\xi$ the
  natural projections.

  Let $\Omega_{M_\xi/k(\xi)}$, $\Omega_{X_\xi/k(\xi)}$ be the sheaves
  of modules of K\"ahler differentials on the generic fibres
  $M_{\xi} = M \otimes_k k(\xi)$ and $X_{\xi} = X \otimes_k k(\xi)$,
  respectively. Since $\iota: X_\xi \to M_\xi^2$ is a closed embedding the
  pull-back 
  $$
  \iota^*: p_1^* \Omega_{M_\xi/k(\xi)} \oplus p_2^*
  \Omega_{M_\xi/k(\xi)} \to \Omega_{X_\xi/k(\xi)}
  $$
  is surjective. Taking stalks at the generic point $\chi$ of $X_\xi$
  we get a surjective map of vector spaces over the field
  $k(\chi) = k(X)$
  $$
  \iota^*: p_1^* \Omega_{M_\xi/k(\xi)} \otimes k(\chi) \oplus p_2^*
  \Omega_{M_\xi/k(\xi)} \otimes k(\chi) \to \Omega_{X_\xi/k(\xi)}
  \otimes k(\chi).
  $$
  Each summand on the left is either trivial or one-dimensional. Since
  $i^*$ is surjective, it follows that at least one of the summands is
  mapped surjectively on the destination. Therefore, the stalk at
  $k(\chi)$ of either
  $\Omega_{X_\xi/k(\xi)}/p_1^* \Omega_{M_\xi/k(\xi)}$ or
  $\Omega_{X_\xi/k(\xi)}/p_2^* \Omega_{M_\xi/k(\xi)}$ vanishes, and we
  conclude.
\end{proof}

Suppose we have a family of pure-dimensional curves
$X \subset M^2 \times T$ such that for some $a \in X_t$ for all
$t \in T$, and assume that for all $t$ in the morphism
$p_1: X_t \to M$ is \'etale in some neighbourhood of $a$. Then by
Lemma~\ref{etale-branch} there exists a unique branch $\alpha$ of $X$
at $a$. It might be the case, though, that $\tau_n(X_t, \alpha)$
vanishes for all $n$, for all $t \in T$, if $p_2$ is everywhere
ramified on the component of $X_t$ that contains $a$. Below we show
that in this case one can consider the family $X\circ X^{-1}$ which
does not have this pathology, and $p_1, p_2$ are both generically
unramified for any of its members.

Recall that if $f: X \to Y$ is a morphism of schemes over a field of
characteristic $p$ then
$\Fr_{f}: X \to X^{(p/Y)}= X \times_{f,Y,\Fr_Y}Y$, the \emph{relative
  Frobenius morphism}, is defined to be $\Fr_X \times f$ where
$\Fr_X, \Fr_Y$ are the absolute Frobenius endomorphisms of $X,Y$,
respectively. If $Y$ is the spectrum of a field then $X^{(p/Y)}$ is
denoted just $X^{(p)}$. If
$X=\Spec R, Y=\Spec S, S = R[r_1, \ldots, r_n]/I$ then
$X^{(p/Y)}= R[r'_1, \ldots, r'_n]/I^{(p)}$ where
$I^{(p)}=\suchthat{ f^{(p)}=\sum_J a_J^p (r')^J }{ f=\sum_J a_J r^J
  \in I}$ (where $J$ is a multiindex), and $\Fr_{X/Y}^*(r'_i) =
r_i^p$. On the level of points, if $X \hookrightarrow Y \times \A^n$
then
$$
\Fr_{X/Y}(y, x_1, \ldots, x_n) = (y, x_1^p, \ldots, x_n^p).
$$
The natural projection $\Fr_{X/Y}(X) \to Y$ is given by
$(y, x_1, \ldots, x_n) \mapsto y$.

\begin{lemma}
  \label{cancel-insep}
  Let $f: X \to Y$ be a finite morphism of irreducible varieties over
  a field of characteristic $p > 0$ and let $F=\Fr_f$ be the relative
  Frobenius morphism.  Assume that $f$ is everywhere ramified. Then
  there exists an $n > 0$ such that the natural projection
  $F^n(X)=X \times_{f,Y,F^n} Y \to Y$ is generically unramified.
\end{lemma}

\begin{proof}
  Since $f$ is everywhere ramified, the field extension
  $k(Y) \subset k(X)$ is inseparable. Let $L$ be the separable closure
  of $k(Y)$ in $k(X)$, then $k(Y) \subset L$ is a separable extension
  and $L \subset k(X)$ is a purely insepearable extension. Since
  $L \subset k(X)$ is a finite extension, there exists a smallest
  number $n$ such that $h^{p^n} \in L$ for any $h \in k(X)$. We claim
  that $k(F^n(X)) \subset L$, which will conclude the proof, as this
  shows that $k(F^n(X))$ is a separable extension of $k(Y)$.
	
  To prove the above claim, let $X_0 \subset X$, $Y_0 \subset Y$ be
  dense open affine subvarieties such that $X_0$ is finite over
  $Y_0$. Then $k[X_0]=k[Y_0][h_1, \ldots, h_n]/I$ and
  $k[F^n(X_0)]=k[Y_0][g_1, \ldots, g_n]/I^{(p^n)}$, and there is an
  embedding of rings $k[X_0] \subset k(X_0)$. It is immediate from the
  definition of the relative Frobenius morphism that there exists an
  injection $k[F^n(X_0)] \hookrightarrow L$ sending $g_i$ to
  $h_i^{p^n}$, so $F^n(X_0)$ is unramified over $Y_0$ and we conclude.
\end{proof}

%

\begin{lemma}
  \label{fr-cancel}
  Let $X \subset M^2 \times T$, $Y \subset M^2 \times S$ be two
  families of pure-dimensional curves. Let us denote projections of
  $M \times M \times T$, resp.  $M \times M \times S$, on products of
  factors by $q$, resp. $q'$, with subscripts. Let $m > 1$ be an
  integer and let $X'=F_{q_{23}}^m(X)$, $Y'=F_{q'_{23}}^m(Y)$. Then
  $$
  X \circ Y^{-1} = X' \circ (Y')^{-1}.
  $$
\end{lemma}

\begin{proof}
  Let us denote projections from
  $M \times M \times M \times T \times S$,
  $M \times M^{(p^m)} \times M \times T \times S$, onto products of
  factors by $r$, resp. $r'$, with subscripts.  After unravelling the
  definitions one observes that
  $$
  X \circ Y^{-1} = r_{1345}(Z) \qquad X' \circ Y'^{-1} =
  r'_{1345}(\Fr_{r_{1345}}(Z))
  $$
  for
  $Z=r_{1245}^{-1}(X \times S) \cap r_{2345}^{-1}(Y^{-1} \times T)
  \subset M^3 \times T \times S$. These projections coincide, since by
  the definition of the relative Frobenius morphism
  $r'_{1345} \circ \Fr_{p_{1345}} = r_{1345}$.
\end{proof}

For the benefit of the reader, let us consider the situation in 
Lemma~\ref{fr-cancel} at the level of points. Denote by
$F: M \to M^{(p)}$ the Frobenius morphism and assume $M$ is affine and
cut out by the equation $f(x_1, \ldots, x_n)$ in $\A^n$, then $M^{(p)}$ is
cut out by $f(x_1^p, ..., x_n^p) = f^p$, and a point $(x_1, ..., x_n)$
is sent by $F$ to $(x_1^p, \ldots, x_n^p)$. The map $\Fr_{q_{23}}$ in
Lemma~\ref{fr-cancel} sends a tuple $(x,y,t) \in M^2 \times T$ to
$(F(x), y, t)$ and similarly for $\Fr_{q'_{23}}$. By definition
$$
\begin{array}{rcl}
  (b,a,t) \in Y \textrm{ if and only if } (F(b), a, t) \in Y', \\
  (b,c,s) \in X \textrm{ if and only if } (F(b), c, s) \in X'. \\
\end{array}
$$

Consider
$$
Z =  \suchthat{(a,b,c,t,s)}{ (b,a,t) \in X, (b,c,s) \in Y}
$$
then
$$
X \circ Y^{-1} = \suchthat{(a,c,t,s) \in M^2 \times T \times S}{
  \exists b\, (b,a,t) \in Y, (b,c,s) \in X = p_{1345}(Z)}.
$$
Also, $\Fr_{r_{1345}}(a,b,c,t,s) = (a,F(b),c,t,s)$ and
$$
\begin{array}{lllll}
X' \circ Y'^{-1} & = & \{ (a,c,t,s) \in M^2 \times T \times
                       S & \mid & \exists b\ (F(b),a,t) \in Y',  \\
  &&&& (F(b),c,s) \in X =  p_{1345}(\Fr_{r_{1345}}(Z)) \}.  \\
\end{array}
$$

\subsection{Interpretation of a one-dimensional group}
\label{int-group}

In the present section we construct a group interpretable in $\Mm$. As
already explained, this will be done by constructing a group
configuration in $\Mm$. In order to construct this group configuration
a one-dimensional algebraic group (Lemma~\ref{many-slopes}) $G$
associated with slopes is `lifted', using
Proposition~\ref{tangency-intersection}, to a group configuration in
$\Mm$. \\

\textbf{Throughout this section and until the end of this paper we fix
  an algebraic curve $M$ over an algebraically closed field $K$ of
  infinite transcendence degree, and a reduct $\Mm$ of the full
  Zariski structure $\M$ on $M$. We assume that the reduct is not
  locally modular. By default the term \emph{definable} will refer to
  definability in $\Mm$. Unless explicitly stated otherwise, by
  definable families we mean \emph{stationary} nearly faithful ample
  families of curves. Where a family $X\to T$ is \emph{stationary} if
  every
  definable open subset of $T$ is dense. }\\


Before we proceed, we need a couple of easy observations: 

\begin{lemma}
  \label{automorphisms} 
  Let $r: \End(k[\eps]/(\eps^{n+1}) \to \End(k[\eps]/(\eps^2)$ be the
  map sending an endomorphism $\varphi$ to the endomorphism
  $\eps\mapsto \varphi(\eps) \mod \eps^2$. Then
  $$
  \Aut( k[\eps]/(\eps^{n+1})) = r^{-1}(\Aut( k[\eps]/(\eps^2))).
  $$
\end{lemma}

\begin{proof}
  Straightforward (see a similar statement for formal power series,
  for example, in \cite[Corollary~7.17]{eisenbud}).
\end{proof}

\begin{lemma}\label{near faithful}
  Let $X\to T$ and $Y\to S$ be one-dimensional definable families of
  strongly minimal subsets of $M^2$. Assume that for all $t\in T$,
  $s\in S$ all projections $p_i: X_t\to M$, $p_i:Y_s\to M$ are
  dominant. Then either for generic $(t,s)$ the set
  $\{(t',s'): |X_{t'}\circ Y_{s'}\cap X_t\circ Y_s|=\infty\}$ is
  finite or $\Mm$ interprets a one-dimensional group.
\end{lemma}

\begin{proof}
  Fix $(t,s)$ generic. Since any curve in $X\circ Y$ intersecting
  $X_t\circ Y_s$ in an infinite set must contain (up to a finite set)
  a strongly minimal component of $X_t\circ Y_s$, and since only
  finitely many such components exists, it will suffice to show that
  any such component is contained in finitely many members of
  $X\circ Y$.
	
  Let $E\subseteq X_t\circ Y_s$ be strongly minimal. By (the proof of)
  \cite[Lemma 3.20]{ElHaPe} either $\Mm$ interprets a one-dimensional
  group or $\dim_\Mm(\mathrm{Cb}_\Mm(E/\0))=2$ (the latter notation can be
  interpreted, equivalently, as: there exist an $\Mm$-definable nearly
  faithful family of curves defined over a two-dimensional parameter set
  and $E$ is generic in that family). We may assume the latter case
  occurs. So, by obvious dimension considerations
  $s,t\in \acl_\Mm(\mathrm{Cb}_\Mm(E/\0))$. So there are only finitely many $(t',s')$
  such that $E\subseteq X_{t'}\circ Y_{s'}$, which is what we had to
  show.
\end{proof}

\begin{rem}
	Recall that our aim in this section is to interpret in $\Mm$ a strongly minimal group $G$. It follows from the previous lemma that one way of achieving this is to find $X\to T$ and $Y\to S$  one-dimensional definable families of strongly minimal subsets of $M^2$ with the property that  $\{(t',s'): |X_{t'}\circ Y_{s'}\cap X_t\circ Y_s|=\infty\}$ is infinite. In order not to overload the formulation of the sequel we will tacitly assume that, whenever Lemma \ref{near faithful} is invoked, this is not the case -- as otherwise we have found our group, and we can move on to the next section. \\
	
\end{rem}

We now proceed to finding the 1-dimensional algebraic group of slopes needed for the construction of the group configuration: 

\begin{lemma}
  \label{many-slopes} 

	
  There exists a nearly faithful definable family
  $Y \subset M^2 \times S$ with $S$ strongly minimal, a locally closed
  irreducible set $S_0 \subset S$, a point $a=(a_1,a_2) \in M^2$,
  $a_1=a_2$, such that $a \in Y_s$ for all $s \in S_0$, and a family
  of branches $\beta$ of $Y \times_S S_0$ at $a$ such that for some
  $n > 0$ the locally closed set
  $$
  \suchthat{ \tau_n(Y_s, \beta_s) } { s \in S_0 }
  $$
  almost coincides with a one-dimensional connected subgroup
  $H \subset \Aut(k[x]/(x^{n+1}))$. 
\end{lemma}

\begin{proof}
  Fix some nearly faithful definable family $X\subseteq M^2\times T$
  witnessing non local modularity of $\Mm$ and such that $X_t$
  is strongly minimal for generic $t\in T$, as provided by Fact
  ~\ref{F:nearlyfaithful}. We may further require that the fibres
  $\pi^{-1}_i(a)$ for both projections of $X_t$ on the factors $M$ are
  finite for all $a\in M$.
  
  Pick an irreducible component $X'$ of $X$ dominant over an
  irreducible component $T_0 \subset T$ of maximal dimension, and such
  that $X'$ is a family of curves. Let $M_0$ be the connected
  component of $M$ such that $M_0^2 \times T_0$ contains $X'$. By
  Lemma~\ref{projections} applied to the closure of $X'$, without loss
  of generality, we may assume that the restriction of $p_1$ to $X'_t$
  is dominant and generically \'etale for $t$ in a dense subset
  $T_1 \subset T_0$. By Lemma~\ref{cancel-insep} there exists a number
  $m$ such that the restriction of $p_{23}$ to
  $X''=\Fr_{p_{23}}^m(X') \cap M_0^2 \times T_1$ is generically
  unramified, and since $X'$ is nearly faithful, the projection is also dominant.
  In particular for any $t \in T_1$ the projection
  $p_2: X''_t \to M_0$ is generically unramified.

  For each $a \in M_0^2$ consider the set $S^a \subset T_1$ of
  $t\in T_1$ such that $a \in X''_t$ and denote
  $X^a=X \cap M^2 \times S^a$. Let $U \subset X''$ be the complement
  of the ramification locus of the restriction of $p_{23}$ to $X''$.
  It follows from dimension considerations that there exists
  $a \in M_0^2$, and an irreducible locally closed subset
  $S_0 \subset S^a$ such that $\dim S_0=1$, $\{a\} \times S_0 \cap U$
  is dense in $\{a\} \times S_0$, and $a \in X''_t$ is smooth for
  $t \in S_0$.  Because $a \in X''_t$ is smooth for any $t \in S_0$,
  there exists by Lemma~\ref{etale-branch} a unique family of branches
  $\alpha$ of $X'' \cap M_0^2 \times S_0$ at $a$.  Then
  $\tau_1(X''_t, \alpha_t) \neq 0$ for $t$ in a dense open subset of
  $S_0$ by the choice of $S_0$, and so by Lemma~\ref{automorphisms}
  $\tau_n(X''_t, \alpha_t) \in \Aut(k[x]/(x^{n+1}))$ for all
  $n \geq 1$ for all $t\in S_0$.  Pick some $t_0 \in S_0$ generic over all the data and let
  $Y = X^a \circ X_{t_0}^{-1}$. Then by Lemma~\ref{fr-cancel}
  $X \circ X_{t_0}^{-1} \cap M_0^2 \times S_0 = X'' \circ
  (X''_{t_0})^{-1}$ and
  $\tau_1(Y_t, \alpha_t \circ \alpha_{t_0}^{-1}) = \tau_1(X''_t \circ
  (X''_{t_0})^{-1}, \alpha_t \circ \alpha_{t_0}^{-1}) \in
  \Aut(k[x]/(x^{n+1}))$ for $t$ in a dense open subset of
  $S_0$. Clearly, $\alpha \circ \alpha_{t_0}^{-1}$ is a family of
  branches at a point $(a_1, a_2) \in M_0^2$ such that $a_1=a_2$.
  	
  By Krull's Intersection theorem and since $S_0$ has non-zero
  dimension, there exists a smallest number $n$ such that
  $|\{\tau_n(X_t, \alpha_t):t\in S_0\}| > 1$. If $n=1$ then
  $\{\tau_n(X_t, \alpha_t):t\in S_0\}$ coincides with a
  one-dimensional subgroup of $\Aut(k[x]/(x^2)) \cong k^\times$ up to a
  finite set.  If $n>1$ then the slope $\tau_{n-1}(X'_t, \alpha_t)$ as
  $t$ ranges in $S_0$ is constant, and therefore
  $\tau_{n-1}(Y_t, \alpha_t \circ \alpha_{t_0}^{-1})=1$. It follows
  that $\tau_n(Y_t, \alpha_t \circ \alpha_{t_0}^{-1})$ almost
  coincides with $\Ker (\Aut(k[x]/(x^{n+1})) \to
  \Aut(k[x]/(x^n)))$. In either case the family $Y$ satisfies the main
  part of the lemma over the irreducible component $S_0$. Near
  faithfulness of $Y$ follows from Lemma~\ref{near faithful} applied
  to $X^a$ and $(X^a)^{-1}$, observing that $Y$ is a subfamily of
  $X^a\circ (X^a)^{-1}$, and that a generic member of $Y$ is generic
  (over $a$, not over $a,t_0$) in $X^a\circ (X^a)^{-1}$.
\end{proof}

Before proceeding to the construction of a group in $\Mm$ we need some
more preliminary work. First, we fix some \textit{ad hoc} terminology and notation that will simplify the discussion: \\

\noindent\textbf{Notation} Let $X\to T$ be a definable family of curves in $M^2$. We denote: 
\begin{enumerate}
	\item For $a\in M^2$ let $T^a:=\{X_t:a\in X_t\}$ be the definable sub-family of all curves incident to the point $a$. 
	\item $X^0:=\{X_t^0: t\in T\}$ where $X_t^0$ is the set of algebro-geometric 0-dimensional components of $X_t$. 
	\item $X^1\subseteq X\times_T T'$ is a family of pure-dimensional curves for a denes $T'\subseteq T$, as provided by Lemma \ref{pure-dim}. 
\end{enumerate}

\begin{defn}\label{def-useful}
  We say that a nearly faithful $\Mm$-definable family of curves
  $X\to T$ satisfies \emph{property $(a,n)$} for $a \in M^2$ and a
  positive integer $n$ if:
  \begin{enumerate}
  \item $a\in X_t$ for all $t$.
  \item There exists a family $\beta$ of branches of $X$ at $a$ such
    that $\{\tau_n(X_t, \beta_t): t\in T\}$ is one-dimensional and
    contains, up to a finite set, a one-dimensional connected
    algebraic group, $H$.
  \item For all $a'\in M^2$, if $a'\neq a$ then $\dim(T^{a'})=0$. 
  \item If $p$ belongs to a zero-dimensional component of
    $X_s\circ X_t$ for $s,t\in T$ $\M$-independent generics then $p$
    is $\M$-generic in $M^2$.
  \end{enumerate}
  The group $H$ is \emph{the group of slopes of $X$ at $a$}
  (associated with the family of branches $\beta$).  
\end{defn}

To show the existence of families that satisfy property $(a,n)$ for
some $a,n$, we need to show:

\begin{lemma}\label{isolated}
  Let $X\to T$, $Y\to S$ be stationary families.
  Then there exists $X'\to T$, $Y'\to S$ $\Mm$-definable over
  $\acl_\M(\0)$ such that
  \begin{enumerate}
  \item $X_t=X'_t$, $Y_s=Y_s'$ up to a finite set, for all $t\in T$, $s\in S$. 
  \item If $t\in T$, $s\in S$ are $\M$-independent
    generics and $a\in X_t\circ Y_s$ is a
    zero-dimensional component then $a$ is $\M$-generic
    over $\0$.
  \end{enumerate}
\end{lemma}
\begin{proof}
  We may assume that for $t\in T$ generic, if $a$ is an isolated
  component of $X_t$ then $a\notin\acl_\M(\0)$. Otherwise, note that
  by stationarity (and genericity of $t$) we get that $a\in X_{t'}$
  for all generic $t'\in T$. So $a\in \acl_\M(\0)$. Since there are at
  most finitely many $b$ incident to all but finitely many $X_{t'}$,
  we may simply set $X':=(X\setminus \{a\})\times T$, eliminating the
  problem in finitely many similar steps.
	
  Similarly, we may assume that if $a=(a_1,a_2)$ is a zero-dimensional component of $X_t$ then $a_1,a_2\notin \acl_\M(\0)$. Thus, we may assume that both $(a_1,a_2)$ are $\M$-generic over $\0$. The same is, of course, true of $Y$. 	
	
  Denoting $X_t^0$, $Y_s^0$ the zero-dimensional components of
  $X_t$, $Y_s$ and noting that
  $(X_t\circ Y_s)^0\subseteq X_t^0\circ Y_s \cup X^t\circ Y_s^0$
  we get that for $s,t\in T$ independent generics any isolated
  point of $X_t\circ Y_s$ is generic over $\0$.
\end{proof}

We have thus shown:

\begin{cor}\label{useful}
  There exists $a_1\in M$, a natural number $n>0$ and a one-dimensional
  definable family of curves that satisfies property $(a,n)$ for
  $a=(a_1, a_1)$.
\end{cor}
\begin{proof}
  Clauses (1) and (2) of the definition of property $(a,n)$ are
  achieved by taking a family as provided by
  Lemma~\ref{many-slopes}. Condition (4) is provided by
  Lemma~\ref{near faithful}, and condition (3) is obtained by
  removing finitely many points common to all generic independent
  curves in the resulting family.
\end{proof}

The same proofs give also: 
\begin{cor}\label{cor-useful}
  If $X\to T$ is a family that satisfies property $(a,n)$, then up to
  -- possibly -- finitely many corrections, $X\circ X$ and
  $X\circ X^{-1}$ also satisfy property $(a,n)$.
\end{cor}

Note however that in the above corollary if $X$ is one-dimensional then
the families $X\circ X$ and $X\circ X^{-1}$ will not be
one-dimensional. It follows, however, that if $t\in T$ is generic then
the one-dimensional families $X\circ X_t$ and $X\circ X_t^{-1}$ will
satisfy property $(a,n)$. The following is a strengthening of the
above observation that we will need later on for technical reasons:

\begin{lemma}\label{near-faithful1}
  Let $X\to T$ be a family that satisfies property $(a,n)$. Let $H$ be
  the group of slopes of $X$ at $a$ (associated with some family of
  branches).  Then there exists a one-dimensional nearly faithful
  family of strongly minimal sets $Z\to L$ such that $a\in Z_l$ for
  all $l$, and there exists a family of branches $\gamma$ at $Q$ such
  that $\tau_n(Z_l,\beta_l)=1\in H$ for all $l \in L$.
\end{lemma}

\begin{proof}
  There exists an $\M$-irreducible component $W\subseteq T$ such that
  $\tau_n(X_t, \beta_t)\in H$ for all $t\in W$ (and in particular the
  slope is defined). Let $t_0\in W$ be generic. So there exists some
  $t_1\in T$ such that
  $\tau_n(X_{t_1},\beta_1)=\tau_n(X_{t_0}, \beta_t)^{-1}$. Let
  $Z_{l_0}\subseteq X_{t_1}\circ X_{t_0}$ be the $\Mm$-definable,
  strongly minimal component containing the branch
  $\beta_{t_1}\circ \beta_{t_0}$. Let $Z\to L$ be the $\Mm$-definable family
  whose generic member is $Z_l$. So there is an $\M$-generic
  sub-family of $Z\to L$ with the property that
  $\tau_n(Z_{l'},\gamma_{l'})=1$ (see
  Proposition~\ref{compo-slope}) for a family $\gamma$ of branches of
  $Z$ at $a$ and for all $l'$ in that sub-family.
\end{proof}

We are finally ready to prove the main result of this section:

\begin{thm}
  \label{first-group} 
  \label{pure-dim-first-group} Let $\Mm$ be a non locally modular
  reduct of an algebraic curve $M$ over an algebraically closed
  field. Then $\Mm$ interprets a one-dimensional group.
 \end{thm} 
 
 \begin{proof}

   We prove the theorem by constructing a group configuration. By
   Corollary~\ref{reduction} we may assume that $M$ is smooth and we
   identify $M$ with $M(K)$ for some algebraically closed field $K$ of
   infinite transcendence degree. We will freely use the remark after
   Definition~\ref{branches} and Lemma~\ref{pure-dim}, referring to
   branches of suitable pure-dimensional subfamilies of definable
   families of curves when we speak about branches of definable
   families of curves.
   
   Let $X\to T$ be a one-dimensional definable family of curves that
   satisfies property $(a,n)$ for some point $a=(a_1,a_1)$ as provided
   by Corollary~\ref{useful}, $H$ the associated group of slopes for
   the family of branches $\beta$. Absorbing into the language the
   parameters needed to define $X$, we may assume that it is
   $\0$-definable.
 	
   We fix a standard group configuration 
   $$
   \mathcal H:=\{   g, h, k, gh, gk, h^{-1}k  \}
   $$
   associated with the action of $H$ on itself by multiplication.
 	
   By Lemma~\ref{many-slopes} there exists an irreducible
   component $W\subseteq T$ such that $\tau_n(X_t, \beta_t)\in H$ for
   all generic $t\in W$.  Identifying (up to a finite set)
   $\{\tau_n(X_t, \beta_t): t\in W\}$ with elements of
   $H\le \mathrm{Aut}(k[\varepsilon]/(\varepsilon^{n+1}))$ we get that
   $\tau_n(X_t, \beta_t)$ is $\M$-inter-algebraic with $t$. At the
   price of replacing $W$ with a (dense) open subset, we may assume
   that $W$ is smooth.
 	
   Any $\M$-independent points $s,t\in W$ generic over all the data
   are, in particular, generic and independent in the sense of the
   reduct $\Mm$. Let $u\in T$ be such that
   $$
   \tau_n(X_u, \beta_u)=\tau_n(X_s, \beta_s)
   \tau_n(X_t, \beta_t)=\tau_n(X_s\circ X_t, \beta_s\circ \beta_t).
   $$
 	
   Such a $u$ exists, since the relative slopes of $X_t$ and $X_s$ are
   generic in $H$, which is one-dimensional. Since the product of two
   independent generic elements of $H$ is again generic in $H$, we can
   find such a $u$.
 	
   Getting back to our group configuration $\mathcal H$ the above
   construction gives us a subset of $T$,
   \[
     \mathcal T_\mathcal H:=\{t_g, t_h, t_k, t_{gh}, t_{gk}, t_{h^{-1}k}\}	
   \]
   such that for every $s\in \mathcal H$ we have
   $\tau_n(X_{t_s}, \beta_{t_s})=s$. Our goal is to show that
   $\mathcal T_{\mathcal H}$ is a group configuration in the sense of
   $\Mm$.

   We have to verify the three sets of conditions appearing in
   Definition~\ref{def-gconf}. That the elements of
   $\mathcal T_\mathcal H$ are pairwise $\Mm$-independent follows from
   the fact that for all $s\in \mathcal H$ also $s\in \acl_\M(t_s)$
   and the elements of $\mathcal H$ are $\M$-independent. That all
   elements in $\mathcal T_\mathcal H$ have dimension $1$ follows from
   the fact that $T$ is strongly minimal and the elements of
   $\mathcal T_\mathcal H$ are generic in $T$. So it remains only to
   verify the third set of conditions, namely, that every collinear
   triple of elements in the following diagram is $\Mm$-dependent:
   $$
   \groupconf{t_g}{t_h}{t_{gh}}{t_k}{t_{gk}}{\; \; \; \; \; t_{h^{-1}k}}
   $$
   The rest of the proof will be dedicated to that end. Since the
   situation is symmetric, it will suffice to show that if $s,t\in W$
   are generic independent
   $\tau(X_u,\beta_u)=\tau_n(X_s,\beta_s)\tau(X_t,\beta_t)$ then
   $u\in \acl_\Mm(s,t)$. Note that since $W$ is $\M$-strongly minimal,
   $u\in \acl_\M(s,t)$.
 	
   To achieve our goal, we would like to apply
   Proposition~\ref{tangency-intersection} to the family
   $\tilde E\to R$ given by $X\circ X$ and the family $X\to T$, in
   order to show that the curve $X_u$ intersects the curve
   $X_s\circ X_t$ in a smaller than generic number of points. The
   problem is that neither $\tilde E\to R$ nor $X\to T$ can be assumed
   to be pure-dimensional families of curves, which is a crucial
   assumption in the statement of the proposition. To circumvent this
   problem, we will show that $X_u\cap (X_s\circ X_t)$ contains no
   zero-dimensional components of either curve, allowing us to apply
   the proposition with the pure-dimensional $\tilde E^1\to R$ and
   $X^1\to T$ without changing the number of intersection points.
 	
   For technical reasons that will be made clear later on, we need to
   slightly twist the family $\tilde E\to R$ that we are working with.
   Let $Z\to L$ be a nearly faithful family of curves that satisfies
   property $(a,m)$ for $m > n$, let $\gamma$ be a family of branches
   at $a$, all as provided by Lemma~\ref{near-faithful1}. We let
   $E\to R$ be $\tilde E\circ Z_{l_0}$ for some $l_0\in L$
   $\M$-generic and independent over all the data. Note that, by
   Proposition ~\ref{compo-slope} and the choice of the family
   $Z\to L$, we have that
   $\tau_n(X_s\circ X_t, \beta_s\circ \beta_t)=\tau_n(X_s\circ
   X_t\circ Z_{l_0}, \beta_s\circ \beta_t\circ \gamma_{l_0})$ whenever
   both sides of the equations are defined. For the sake of clarity we
   let $\alpha$ be the family of branches of $E$ at $a$. Namely,
   $\alpha =\beta\circ \beta \circ \gamma_{l_0}$.
 	 	
   It will be convenient to already note at this stage the following
   slight strengthening of Lemma~\ref{near faithful}:\\
	
   \noindent\textbf{Claim 1} We may assume that if $r\in R$ is generic
   then $|\{r': \tilde E_r\cap \tilde E_{r'}\}|=\infty$ is finite.  

   \proof The claim would follow from Lemma~\ref{near faithful}, if the
   members of $X$ were strongly minimal. In the general case, if $r\in R$ is
   generic and $\tilde E_r=X_s\circ X_t$ then any strongly minimal 
   $F_r\subseteq \tilde E_r$ is contained in $C_s\circ D_t$ for some
   strongly minimal $C_s\subseteq X_s$ and $D_t\subseteq X_t$. By
   Lemma~\ref{near faithful} applied to the families $\{D_t: t\in T\}$
   and $\{C_s: s\in S\}$, we get that $s,t\in \acl_{\Mm}(\Cb(F_r)])$. 
   Since $\Cb(F_r)\in \acl(\Cb(E_r))$ we conclude that $s,t\in \acl(\Cb(E_r))$, which is what we needed. \qed$_{\text{Claim }1}$. \\
	
   Note that the fact that $\tilde E$ is the composition of two copies of $X$ did not play any role in the proof above, and we could invoke Lemma \ref{near faithful} with $X$ and $(X\circ Z_{l_0})$ to get the same conclusion for the family $E:=X\circ (X\circ Z_{l_0})$. \\

   Let us fix some additional notation. We let $R(u)$ be the set of
   all $r\in R$ such that $E_{r}$ is tangent to $X_u$ at $a$, i.e.,
   $\tau_n(E_r, \alpha_r)=\tau_n(X_u,\beta_u)$.  Let $E(u):=\{E_{r}: r\in R(u)\}$. In other words,
   $R(u)$ is the parameter set of all curves in the family $E\to R$
   $n$-tangent to $X_s\circ X_t$ at $a$ and $E(u)$ is the subfamily of $E$
   over the parameter set $R(u)$. So $E(u)\to R(u)$ is an
   $\M$-definable subfamily of $E$ of dimension 1. Fix, once and for all,
   $r=(s,t)\in R$. So $r\in R(u)$ and it is $\M$-generic as
   such. 
   Replacing, if needed, $R(u)$ with the $\M$-definable strongly minimal component of $R(u)$ containing $r$ we may assume that $R(u)$ is strongly minimal. \\
	

	The following is the main step in the proof: \\
	
	\noindent\textbf{Claim 2:} Assume that $u\notin\acl_\Mm(s,t)$.  Let $\{x_1,\dots, x_k\}=(X_u\cap E_r)\setminus \{a\}$. Then 
	$x_i$ is $\M$-generic in $X_u$ for all $i$. 
	
	\proof First, note that $x_i\notin \acl_{\Mm}(\0)$, because
        otherwise, since $u\in T$ is generic, we would get that
        $\dim(T^{x_i})=1$, contradicting clause (3) of
        Definition~\ref{def-useful}. Note that the exact same argument
        shows that $x_i\notin \acl_\M(\0)$. Next, as
        $r\in \acl_\Mm(s,t)$ and since $X$ is an $\Mm$-strongly
        minimal family, our negation assumption implies that $u$ is
        $\Mm$-generic over $r$ and by Lemma~\ref{near faithful} (and the remark following it) applied to the $\acl(\0)$-definable strongly minimal subsets of $E_r$ we get that
        $\dim(\Cb(E_r)/\0)=\dim(r/\0)=2$.

	Now assume that $x_1$ is not $\M$-generic in $X_u$. Since $r$
        is $\M$-generic in $R(u)$ (and thus also in $R$), it follows
        that $\dim(R(u)^{x_1})=1$. Indeed, since $\dim_\Mm(x_1/u)=0$
        (by assumption), it follows that
        $\dim_\M(r/ux_1)=\dim_\M(r/u)$, so $r$ is generic in $R(u)$
        over $x_1$, and the strong minimality of $R(u)$ 
	implies that $x_1\in E_{r'}$ for all generic $r'\in
        R(u)$. Thus, in fact $R(u)$ is a generic subset of $R^{x_1}$.  Recall, moreover, that there exists a family $\alpha$ of branches of all curves in $E(u)$ at $a$ such that $\tau_n(E(u), \alpha_r)=\tau(X_u,\beta_u)$. 

	We will show that this gives the desired conclusion. We split
        the argument into cases according to $\dim_\Mm(x_1/\0)$. The
        case $x_1\in \acl_{\Mm}(\0)$ has already been discarded. If
        $x_1$ is non-$\Mm$-generic in $M^2$ then there exists a curve
        $F$, $\Mm$-definable over $\0$ such that  $x_1$ is generic in
        $F$. So $u$ is contained in the set of all $u'$ such that
        $F\cap X_{u'}\cap E_r\neq \0$. Because $|E_r\cap F|<\infty$
        and  condition (3) of Definition~\ref{def-useful}, there are only finitely many such $u'$. So $u$ is $\Mm$-algebraic over $r$ -- contradicting our assumption. 

	Thus, we may assume that $x_1$ is $\Mm$-generic in $M^2$. We
        will now focus on the family $E\to R$. Since $x_1$ is
        $\Mm$-generic in $M^2$, for any $r_1,r_2\in R^{x_1}$
        independent generic, $m:=|E_{r_1}\cap E_{r_2}|$ is obtained on
        an $\M$-generic subset of parameters of $R\times R$. Consider the
        $\M$-definable family $E^1\to R$ of pure-dimensional curves
        associated with $E\to R$. Note that for $\M$-generic
        independent $u,w\in R$ we have $|E_u\cap E_w|=|E^1_u\cap
        E^1_w|=m$.

        On the other hand, Lemma~\ref{int-flat}, and hence
        Proposition~\ref{tangency-intersection}, is applicable to two
        copies of the family $E^1 \to R$, possibly after shrinking $R$
        so as to ensure, using Fact~\ref{generic-flatness}, that
        $E^1 \to R$ is flat and that $R$ is smooth.

        So, by Lemma~\ref{int-flat} and
        Proposition~\ref{tangency-intersection} there is a dense open
        set $W_0\subseteq R$ defined over $\0$ such that (keeping the above
        notation) for all $(v,w) \in W:=W_0\times W_0$, if
        $\tau_n(E^1_v, \beta_v) = \tau_n(E^1_w,\beta_w)$, then either
        $\dim(E^1_w\cap E^1_v)=1$ or $\#(E^1_w \cap E^1_v ) < m$. 
        
     	For generic $v$ the set of $w$ such that  $(v,w)\in W$ and $\dim(E^1_w\cap E^1_v)=1$ is finite by Claim 1. 
     	So, for generic $(v,w)\in W$ we see that $\dim(E^1_w\cap E^1_v)=0$ so necessarily  $\#(E^1_w \cap E^1_v) < m$. We now show that this must imply that $E^0_w\cap E_r\neq \0$ for generic $w$. Indeed, since $W_0$ is dense in $R$ and $\0$-definable, $\M$-genericity of $r$ in $R$ implies that it is also generic in $W_0$. Since $R(u)^{x_1}$ is
        generic in $R^{x_1}$ (in the sense that it contains an open
        subset of $R^{x_1}$) we can find some $w\in R(u)\cap W_0$
        $\M$-generic and $\M$-independent from $r$ (over all the data gathered so far) so that $(r,w)$ is $\M$-generic in $W$. Moreover,  by definition of
        $R(u)$ we know that
        $\tau_n(E^1_r, \beta_r) = \tau_n(E_w^1,\beta_w)$, and by what we
        have just said, this must imply that
        $\#(E^1_w \cap E^1_r ) < m$. Because $x_1$ is $\Mm$-generic in
        $M^2$ and $w,r\in R^{x_1}$ are $\Mm$-independent generics, they
        are, in fact, independent generic in $R^2$ over $\0$. So
        $\#(E_w \cap E_r ) =m$, implying that $E_w\cap E_r^0\neq
        \0$. 
        
        Finally, since $w$ was $\M$-independent from $r$ and $\M$-generic
        in $R(u)$, and since $E_r^0\subseteq \acl_\M(r)$, we get --
        precisely as above -- that there is some $c\in E_r^0$ such
        that $R^c$ contains $R(u)$, up to a finite set. This implies that $\dim_\M(c/u)=0$, and therefore $\dim_\M(c/\0)\le \dim_\M(u/\0)=1$. This contradicts Corollary~\ref{cor-useful} (specifically, clause (4) of Definition~\ref{def-useful}). \qed$_{\text{Claim }2}$. \\
	
	It follows from Claim 2 that $X_u^0\cap E_r=\0$. We also need
        to show that $X_u$ does not meet $E_r$ in an isolated point of
        the latter. It is here that the twist of the family $\tilde E\to R$ by a generic curve from $Z\to L$ plays its role: \\
	
	\noindent\textbf{Claim 3:} If $u\notin \acl_\Mm(s,t)$ then
        $X_u\cap E_r^0=\0$.  
	\proof Recall that $E_r=X_s\circ X_t\circ Z_{l_0}$. Assume
        that there exists some
        $x_i\in X_u\cap (X_s\circ X_t\circ Z_{l_0})^0$. By
        Lemma~\ref{isolated} applied to $\tilde E(u)\to R(u)$ and
        $Z\to L$, if $r'\in R(u)$ is generic and $l\in L$ is generic
        independent from $r'$, then any $x_i\in (E_{r'}\circ Z_l)^0$ is
        either $\M$-generic over $\0$ or contained in one of finitely
        many sets of the form $\{a\}\times M$ and $M\times \{a\}$ for
        $a\in \acl_\M(0)$.
	But $X_u\cap (M\times \{a\}\cup \{a\}\times M)\subseteq \acl_{\M}(u)$, so $x_i\in \acl_\M(u)$, contradicting the previous claim. \qed$_\text{Claim 3}$\\
	
	The conclusion of the discussion, up to this point, is that if $u\notin \acl_\Mm(s,t)$ then $E_r\cap X_u=E_r^1\cap X_u^1$. This allows us to conclude that, in fact: \\
	
	\noindent\textbf{Claim 4:} $u$ is $\Mm$-algebraic over $t,s$. 
	\proof Assume not. By Proposition~\ref{compo-slope}
        $\tau_n(X_t \circ X_s\circ Z_{l_0}, \beta_t \circ \beta_s\circ
        \gamma_{l_0}) = \tau_n(X_t, \beta_t) \circ \tau_n(X_s,
        \beta_s)$.  Let
        $m = \max_{\bar t, \bar s, \bar u \in T} \# (X_{\bar t} \circ
        X_{\bar s}\circ Z_{l_0} \cap X_{\bar w})$, then by
        Lemma~\ref{sc},
        $m = \#((X_t \circ X_s\circ Z_{l_0}) \cap X_u)$ for $(t,s,u)$
        generic in $T \times T \times T$. Let $\tilde T\subseteq T$ be
        as provided by Lemma~\ref{many-slopes}. By
        Lemma~\ref{int-flat} and
        Proposition~\ref{tangency-intersection}, the set of parameters
        $w \in \tilde T$ such that
        $\tau_n(X_w^1, \beta_w) = \tau_n(X_t^1, \beta_t)\circ
        \tau_n(X_s^1, \alpha_s)$ is contained in
	$$
	W_1:=\suchthat{ w \in T }{ \dim(X_t \circ X_s \cap X_w)=1
          \textrm{ or }\#((X_t \circ X_s\circ Z_{l_0})^1 \cap X^1_w )
          < m }.
	$$
	By strong minimality of $T$ the set $\{w: \#(X_t \circ
        X_s\circ Z_{l_0} \cap X_w ) < m\}$ is finite. Also, for
        $\M$-generic $w$ we have $(X_t \circ X_s\circ Z_{l_0})^1 \cap
        X_w^1=X_t \circ X_s\circ Z_{l_0} \cap X_w$. So the set of $w$
        such that $\#((X_t \circ X_s\circ Z_{l_0})^1 \cap X^1_w ) < m$
        is finite. Since $X$ satisfies property $(a,n)$, by Lemma~\ref{near faithful} the set  $\{w: \dim(X_t \circ X_s\circ Z_{l_0} \cap X_w)=1 \}$ is finite. So $W_1$ is finite. 
	Similarly, 
	\[
	 W:=\suchthat{ w \in T }{ \dim(X_t \circ X_s\circ Z_{l_0} \cap X_w)=1 \textrm{ or
		}\#(X_t \circ X_s\circ Z_{l_0} \cap X_w ) < m }
	\]
	is finite, and moreover, $W$ is $\Mm$-definable. Our assumption that $u\notin \acl_\Mm(s,t)$ allows us to apply Claim 2 combined with Claim 3 to get that $X_s\circ X_t\circ Z_{l_0} \cap X_u=(X_s\circ X_t\circ Z_{l_0})^1 \cap X_u^1$. Since $u\in W_1$ it follows that $u\in W$, proving that, in fact $u\in \acl(s,t)$.
	  \qed$_\text{Claim 4}$\\
	
	Claim 4 shows that, indeed, $\mathcal H_\mathcal T$ is an $\Mm$-group configuration, and the desired conclusion is obtained by applying Fact~\ref{grconf}.

\end{proof}

The next proposition can be proved in greater generality (and follows,
essentially, from \cite[Section 3]{hr-pillay94}), but we only need the
following elementary result:

\begin{prop}
  \label{ga-not-gm}

  In the notation of the previous proof, assume that the group $H$
  almost coinciding with $\{\tau_n(Y_t, \alpha_t): t\in \tilde S\}$ is
  isomorphic to $\G_a$. Then the connected component of the identity
  of the group from the conclusion of the theorem is not
  $\M$-isomorphic to $\G_m$.
\end{prop}

\begin{proof}
In this proof we will be working solely in $\M$.  It suffices to show that if
  $$
  G_1=\set{a, b, a+b, x, x+a, x+b} \textrm{ and }
  G_2=\set{e, f, e \cdot f, y, e \cdot y, f \cdot y}
  $$
  are group configurations for the groups $\G_a$ and $\G_m$,
  respectively, and
  $\acl(a) = \acl(e), \acl(b) = \acl(f), \acl(x)=\acl(y)$ then $G_1$
  and $G_2$ are not inter-algebraic.

  Indeed, $\dim(a,b,e+f,ef)=\dim(a,a+b, ef)$, since $ef$ is inter-algebraic
  with $b$ over $a$. On the other hand, $\dim(a,a+b,ef)=\dim(a+b,ef)$
  since $a+b$ is inter-algebraic with $a$ over $ef$.
\end{proof}

\subsection{Interpretation of the field and proof of the main theorem}
\label{int-field}

In this section we interpret the field $K$ in the reduct $\Mm$, concluding the proof of the main theorem of this paper. The results of the previous subsection allow us replace $\Mm$ with an algebraic group, $G$, interpretable in $\Mm$ (we only have to verify that the induced structure is non-locally modular). As in the previous subsection, the interpretation of the field boils down to the construction of a field configuration. The construction of the field configuration will depend on whether the (connected component of the) group $G$ is isomorphic (in $K$) to $\mathbb G_a$, $\mathbb G_m$ or to an elliptic curve. The question to address is  how to find an $\Mm$-definable strongly minimal $Z\subseteq G^2$ whose set of slopes $\{\tau_1(Z,z): x\in Z\}$ (see below) is infinite. The easiest is the case of an elliptic curve: 

\begin{lemma}
  \label{slopes-elliptic} 
  Let $E$ be an elliptic curve and $Z$ be a closed one-dimensional
  irreducible subset of $G=E^2$. Identify $T_g G$ with $T_0 G$ via the
  isomorphism $d\lambda_g: T_0G \to T_gG$, for
  $\lambda_g(x)=g \cdot x$.  Suppose that for any $z \in Z$ the
  tangent space $T_z Z \subset T_0G$ is constant. Then $Z$ is a coset
  of a closed subgroup of $G$.
\end{lemma}

\begin{proof}
  Since $Z$ is a projective curve with a trivial tangent bundle, it is
  an elliptic curve itself by the Riemann-Roch formula. Since any
  morphism between Abelian varieties with finite fibres which
  preserves the identity automatically preserves the group structure
  by the Rigidity Lemma (see \cite[p.~43]{mumford}), $Z$ is a coset of
  an Abelian subvariety of $G$.
\end{proof}

Let $M$ be an algebraic curve, and consider a curve $Z \subset
M^2$. For every point $z \in Z$ such that $p_1$ is \'etale in a
neighbourhood of $z$, there exists by Lemma~\ref{etale-branch} a
unique branch at $z$, call it $\alpha_z$. We will use the notation
$\tau_n(Z, z) := \tau_n(Z, \alpha_z)$.  For any group $(G, \cdot)$
with identity $e \in G$, for any $a=(a_1, a_2)\in G^2$ define the maps
$t_a: G^2 \to G^2$
$$
t_a(x_1, x_2) = (a_1^{-1 }\cdot x_1, a_2^{-1}\cdot x_2)
$$
and for any one-dimensional locally closed subset
$Z \subset G^2$ define the set $s_n(Z) \subset k[x]/(x^{n+1})$
$$
s_n(Z) = \suchthat{\tau_n(t_z(Z), (e,e)) }{z \in Z}.
$$
Also for any $c \in k$ define $u_c: \G_a^2 \to \G_a^2$
$$
u_c(x_1, x_2) = (x_1, x_2 - c\cdot x_1).
$$
We can consider the family of translates $\bar Z \subset M^2 \times Z$
such that $\bar Z_a = t_a(Z)$. If $Z$ is a pure-dimensional curve then $\bar Z$ is a family of pure-dimensional curves, and if the
projection on the first coordinate is generically \'etale then for a
dense subset $Z_0 \subset Z$ there is a unique family of branches
$\bar\alpha$ of the family $\bar Z|_{Z_0}$ at $(0,0)$ such that
$\bar\alpha_a=\alpha_a$ for all $a \in Z_0$.

\begin{lemma}
  \label{diff-coset}
  Let $Z \subset \G_a^2$ be an irruducible pure-dimensional curve that
  is not contained in a coset of a subgroup of $\G_a^2$. Then either
  $Z - t_x Z$ is not contained in a subgroup of $\G_a^2$ for some
  $x \in Z$, or the closure of $Z$ is cut out by an equation of the
  form $y - ax^2 = 0$, when the characteristic $p$ of the ground field
  is $0$, or, when  $p>0, p\neq 2$, by an equation
  $$
  y - a x^{2p^n} = 0,
  $$
  for some non-negative integer $n$ and a constant $a \in k$.
\end{lemma}

\begin{proof}
  We may assume that $(0,0) \in Z$, and that the projection on the
  first coordinate in $\G_a^2$ is \'etale in a neighbourhood of
  $(0,0) \in Z$.  Assume that $Z - t_x Z$ is contained in a subgroup
  of $\G_a^2$ for all $x \in Z$, and let the closure of $Z$ be cut out
  by a polynomial equation $f(x,y)$. Let $Z' \subset \G_a^4$ be the
  set cut out by the equation $f(u,v)$, and let $Z'' \subset Z'$ be
  cut out in $Z'$ by the equation $f(x - u, y - v)$ where $x,y,u,v$
  are the coordinates. A fibre of $Z''$ over $(u,v) \in \G_a^2$ is
  thus the closure of $Z - t_{(u,v)} Z$. We identify the completion of
  the local ring of $Z'$ at $(0,0,0,0)$ with $k[[x,y,u]]$ via the
  projection on the first three coordinates. If the local equation of
  $Z$ in $\G_a^2$ at $(0,0)$ is $y - g(x), g \in k[[x]]$ then one
  readily sees that the local equation of $Z''$ in $Z'$ at $(0,0,0,0)$
  is
  $$
  h = y - g(u) - g(x) + g(x - u).
  $$
  If $Z - t_x Z$ is contained in a subgroup for all $x \in Z$ then the
  above expression must be the local equation of the subgroup, that
  is, it is of the form $a y^{p^m} - b x^{p^n}$ for some non-negative
  integers $n,m$. Then one sees immediately that $h$ can only be of the
  form $y - b x^{p^n}$, so $g$ can only be of the form
  $g = a  x^2$, when $p=0$, or $g = a  x^{2p^n}$
  otherwise. In particular, if $p=2$ then $Z$ is a subgroup of
  $\G_a^2$.
\end{proof}

\begin{lemma}
  \label{slopes-ga}
  Let $G$ be an algebraic group such that the connected component of
  the identity $G_0$ is isomorophic to $\G_a$. Let $Z \subset G^2$
  be a curve that is not a Boolean combination of cosets of subgroups
  of $G^2$.  Then there exists a pure-dimensional curve
  $W \subset G^2$ definable in $(G, Z)$, and an irreducible component
  $W' \subset W \cap G_0^2$ of dimension 1, such that
  $\dim s_1(W')=1$. In characteristic 0, one can take $W=Z$.
\end{lemma}

\begin{proof}  
  Replacing $Z$ by a shift, we may assume that there exists a
  pure-dimensional irreducible curve $Z_0 \subset Z$ such that
  $(0,0) \in Z_0$, such that $Z_0$ is not contained in a coset, and
  such that the projection on the
  first coordinate in $\G_a^2$ is \'etale in a neighbourhood of
  $(0,0) \in Z$.We
  will construct a finite sequence of curves $W_i \subset \G_a^2$,
  $W_0=Z_0$ that are irreducible components of curves definable in
  $(G, Z)$, such that
  $$
  N(W_i) := \min \suchthat{ n \geq 1 }{ \dim s_n(W_i) = 1 }
  $$
  is strictly decreasing. 

  We pick the local coordinate systems at all points $c \in \G_a$ in a
  uniform way, with $k[[x]] \cong \Oo_{\G_a, c}$ given by
  $x \mapsto x-c$.  Clearly, $\tau_n(Z, x)=\tau_n(t_b(Z), t_b(x))$ for
  any $b \in \G_a^2$ and $u_c(Z) = Z + L_c$ where $L_c$ is the line
  $x_2 = -c \cdot x_1$. By Proposition~\ref{sum-slope}
  $$
  \tau_1(u_c(Z), (0,0)) = \tau_1(Z, (0,0)) - c.
  $$
  If $N(W_i) > 1$, then, again by Proposition~\ref{sum-slope} 
  $$
  \tau_1(W_i - t_x(W_i), (0,0)) = 0
  $$ for all $x \in W_i$ and since
  $W_i$ was not a coset, by Lemma~\ref{diff-coset}, either
  $W_i - t_x(W_i)$ is not a subgroup for some $x \neq 0$, or it is cut
  out by an equation of the form $y - a x^2$ or $y - ax^{2p^n}$ for
  some non-negative integer $n$, $a \in k$. In the latter case by
  Lemma~\ref{fr-cancel} we have that $N(W_{i+1}) = 1$ for
  $W_{i+1} = W_i \circ (t_x (W_i))^{-1}$ for some $x \in W_i$. In the
  former case, pick an $x$ such that $W_i - t_x(W_i)$ is not a
  subgroup of $\G_a^2$ and define $Y_i=W_i - t_x(W_i)$.  If
  $s_1(W_i) = \{c\}$ then $s_1(Y_i) = \{0\}$ and since $Y_i$ is not a
  coset, $p_2$ restricted to $u_c(Y_i)$ is dominant and everywhere
  ramified. In particular, since the latter is impossible when the
  characteristic of the base field is 0, it follows in that case that
  $N(W_0)=1$.

  By Lemma~\ref{cancel-insep}, there exists a number $m$ such that
  $p_2$ restricted to $Y'_i=\Fr^m_{p_2}(Y_i)$ is generically
  unramified. Let $a \in Y_i$ be a point such that $p_1$ is \'etale
  over $M_0$ in some neighbourhoods of $a$, and let the local
  equation of $Y_i$ at $a, b$ be $x_2 - f^{p^m}(x_1)$,
  respectively. Then the local equation of $Y'_i$ at the point
  $\Fr^m_{p_2}(a)$ is $x_2 - f(x_1)$.
  Define $W_{i+1} = Y_i \circ Y_i^{-1}$, 
  then, by Lemma~\ref{fr-cancel}, $W_{i+1} = Y'_i \circ (Y'_i)^{-1}$ and
  $$
  \tau_n(W_{i+1}, (a_2, b_2)) = \tau_{p^m}(Y_i, a) \circ
  (\tau_{p^m}(Y_i, a))^{-1}
  $$
  for all points $a=(a_1, a_2), b=(b_1, b_2) \in Y_i$ such that
  $a_1=b_1$ and such that the right-hand side makes sense. In
  particular, $p_2$ restricted to $W_{i+1}$ is generically unramified
  and it follows that $N(W_{i+1}) = N(W_i)/p^m < N(W_i)$. Therefore
  for some finite $l$, $N(W_l)=1$. 
\end{proof}

\begin{lemma}
  \label{slopes-gm}  
  Let $G$ be a one-dimensional algebraic group such that the connected
  component of the identity $G_0$ is isomorphic to $\G_m$. Let
  $Z \subset G^2$ be a curve that is not a Boolean combination of
  cosets of subgroups of $G^2$.  Then either there exists an
  irreducible pure-dimensional curve $Z_0 \subset Z$ such that
  $\dim s_1(Z_0)=1$, or there exists a group definable in $(G, Z)$
  such that its connected component of the identity is not isomorphic
  to $\G_m$.
\end{lemma}

\begin{proof}
  Pick a local coordinate systems on $\G_m$, uniformly, as in the
  proof of Lemma \ref{slopes-ga}.  Assume that $\dim s_1(Z)=0$, and so
  $s_1(Z_i)$ is a singleton for each one-dimensional irreducible
  component $Z_i \subset Z$.  Let $Z_0$ be one of the irreducible
  components of $Z$ that is not contained in a coset, then there
  exists a  smallest $n > 1$ such that $\dim s_n(Z_0)=1$. Then, by the
  same reasoning as in the proof of Lemma~\ref{many-slopes}, we may
  consider the family $Y \subset G^2 \times Z$ by putting
  $Y_z = t_z(Z) \circ (t_{z_0}(Z))^{-1}$ for some $z_0 \in Z_0$, so
  that $\{ \tau_1(Y_z, (0,0)) \mid z \in Z_0 \}$ almost coincides with
  $$
  \Ker (\Aut(k[x]/(x^{n+1}) \to \Aut(k[x]/(x^n)) \cong \G_a.
  $$
  The definable family $Y$ can be used to construct a group
  configuration as in the proof of Theorem~\ref{first-group}, and
  therefore a group is interpretable in $(G,Z)$. By
  Proposition~\ref{ga-not-gm}, the connected component of the identity
  of this group is not isomorphic to $\G_m$.
\end{proof}

We can finally interpret the field: 

\begin{thm}
  \label{second-group}

  Let $G$ be a one-dimensional algebraic group over an algebraically
  closed field, $Z \subset G^2$ be a one-dimensional constructible
  subset that is not a Boolean combination of cosets. Then
  $(G,\cdot,Z)$ interprets a field.
\end{thm}

\begin{proof}
  Let $G_0$ be the connected component of the identity $e$ of $G$.  If
  $G_0=\G_a$ or $G_0$ is an elliptic curve then, by
  Lemmas~\ref{slopes-elliptic},~\ref{slopes-ga}, there exists a
  definable family $Y \subset G^2 \times S$ of curves, $S$ strongly
  minimal, and an irreducible locally closed set $S_0 \subset S$ such
  that there is a unique family of branches $\alpha$ of
  $Y_0 = Y \cap S_0$ at $(e,e) \in G^2$, and such that
  $\tau_1(Y_s, \alpha_s)$ is not constant as $s$ ranges in $S_0$. By
  Lemma~\ref{slopes-gm}' either such a family exists, or a definable
  one-dimensional group $G'$ with the connected component of the
  identity not isomorphic to $\G_m$ (and therefore isomorphic to
  either $\G_a$ or to an elliptic curve) is interpretable in
  $(G,\cdot,Z)$, and we may prove the theorem for the structure
  induced on $G'$. We, therefore, may continue with the assumption
  that such a family exists. Clearly, $Y_0 \to S_0$ and we may assume
  that $S_0$ is smooth at the price of possibly shrinking
  $S_0$. Further shrinking $S_0$ we can ensure $Y_0 \to S_0$ to be
  flat (by Fact~\ref{generic-flatness}). Let $K$ be a field of
  infinite transcendence degree over the base field $k$. We identify
  first order slopes, which are truncated polynomials in
  $k[\eps]/(\eps^2)$ divisible by $\eps$, with $K$, and we will use
  multiplicative notation for composition. We will freely use the
  remark after Definition~\ref{branches} and Lemma~\ref{pure-dim},
  referring to branches of suitable pure-dimensional subfamilies of
  definable families of curves when we speak about branches of
  definable families of curves.
  
  Take $a_1,a_2,b_1,b_2, u \in S_0(K)$ generic and pairwise
  independent. Let $c_1,c_2 \in S_0(K)$ be such that
  \begin{eqnarray*}
    \tau_1(Y_{c_1}, \alpha_{c_1}) & = & \tau_1(Y_{a_1}, \alpha_{a_1})  \tau_1(Y_{b_1}, \alpha_{b_1})\\
    \tau_1(Y_{c_2}, \alpha_{c_2}) &                                                                       = & \tau_1(Y_{a_2}, \alpha_{a_2})  \tau_1(Y_{b_1}, \alpha_{b_1}) +  \tau(Y_{b_2}, \alpha_{b_2}).
  \end{eqnarray*}
  This is possible, since the image of the function
  $s \mapsto \tau_1(Y_s, \alpha_s)$ for $s$ ranging in $S_0$ is of
  dimension 1, and the values of slopes on the right hand side of the
  equations above are generic in $\End(k[x]/(x^2))$ for generic
  pairwise independent values of parameters. Therefore
  $$
  \tau_1(Y_{a_1}, \alpha_{a_1}) \tau_1(Y_{b_1}, \alpha_{b_1}) \textrm{
    and } \tau_1(Y_{a_2}, \alpha_{a_2}) \tau_1(Y_{b_1}, \alpha_{b_1})
  + \tau_1(Y_{b_2}, \alpha_{b_2})
  $$
  are generic, and $c_1, c_2$  as required can be found in $S_0(K)$.	
  Let $z, v$ be such that
  \begin{eqnarray*}
    \tau_1(Y_z, \alpha_z) & = & \tau_1(Y_{a_1}, \alpha_{a_1})
                                \tau_1(Y_u, \alpha_u) +
                                \tau_1(Y_{a_2}, \alpha_{a_2}),\\ 
    \tau_1(Y_v, \alpha_v) & = & \tau_1(Y_{b_1}, \alpha_{b_1})^{-1}
                                \tau_1(Y_u, \alpha_u) -
                                \tau_1(Y_{b_2}, \alpha_{b_2}). 
  \end{eqnarray*}
  By a similar reasoning, $z, v$ are generic. It also follows from the
  way $c_1, c_2, z, v$ were defined that
  $$
  \tau_1(Y_z, \alpha_z)  =  \tau_1(Y_{c_1}, \alpha_{c_1}) \tau_1(Y_v,
  \alpha_v) +  \tau_1(Y_{c_2}, \alpha_{c_2}).
  $$
	
  We will now show that $(c_1, c_2)$ is algebraic over $(a_1, a_2)$
  and $(b_1, b_2)$ in the sense of $(G,\cdot,Z)$.  By
  Propositions~\ref{compo-slope} and \ref{sum-slope},
  \begin{eqnarray*}
    \tau_1(Y_{a_1} \circ Y_{b_1}, \alpha_{a_1} \circ \alpha_{b_1}) & =
    & \tau_1(Y_{a_1}, \alpha_{a_1}) \tau_1(Y_{b_1}, \alpha_{b_1}),\\ 
    \tau_1(Y_{a_2} \circ Y_{b_1} + Y_{b_2}, \alpha_{a_2} \circ
    \alpha_{b_1} + \alpha_{b_2}) & = & \tau_1(Y_{a_2}, \alpha_{a_2}) \tau_1(Y_{b_1}, \alpha_{b_1}) + \tau_1(Y_{b_2}, \alpha_{b_2}).
  \end{eqnarray*}
	
  Let
  $l_1 = \max_{c_1, a_1, b_1 \in S_0} \#(Y_{c_1} \cap Y_{a_1} \circ
  Y_{b_1})$, and let
  $l_2 = \#(Y_{c_2} \times_{G^2} Y_{a_2} \circ Y_{b_1} + Y_{b_2})$ for
  $a_1, a_2$, $b_1, b_2$, $c_1, c_2 \in S_0$ generic and
  independent. Since the number of intersection points is an
  $(G,\cdot, Z)$-definable property, it does not matter what
  particular parameters $a_i, b_i, c_i$ we take as long as they are
  generic and independent (in the sense of $(G,\cdot, Z)$).  By
  Lemma~\ref{int-flat} and Proposition~\ref{tangency-intersection}, the
  $(M,X)$-definable set
  $$
  \suchthat{w \in S_0}{ \dim (Y_w  \cap (Y_{a_1} \circ Y_{b_1})) = 1
    \textrm{ or } \#(Y_w  \cap (Y_{a_1} \circ Y_{b_1})) < l_1}
  $$
  contains $c_1$ and by definition of $l_1$ is finite. By
  Lemma~\ref{int-flat} and Proposition~\ref{tangency-intersection}
  again, the $(M,X)$-definable set
  $$
  \suchthat{w \in S_0}{ \dim (Y_w \cap (Y_{a_2} \circ Y_{b_1} +
    Y_{b_2})) = 1 \textrm{ or } \#(Y_w \cap (Y_{a_2} \circ
    Y_{b_1} + Y_{b_2})) < l_2}
  $$
  contains $c_2$ and by definition of $l_2$ is finite.
	
  Arguing in a similar fashion, by application of
  Lemma~\ref{int-flat} and
  Proposition~\ref{tangency-intersection}, we deduce that for all
  lines in the diagram
  $$
  \groupconf{\hspace{-2em}(a_1, a_2)}
  {\hspace{-2em}(b_1, b_2)}
  {\hspace{-2em}(c_1, c_2)}{z}{u}
  {\hspace{1em}v}
  $$
  each vertex is in the algebraic closure of two other collinear
  vertices, and so this constitutes a group configuration. Therefore,
  by Fact~\ref{grconf}, there exists a two-dimensional group definable
  in $(G,\cdot,Z)$ that acts transitively on a one-dimensional set.
  
  The conditions of the Fact~\ref{grconf-faithfulness} are verified as
  well: for instance, for the uppermost line,
  $B=\{ \tau_1(Y_{a_1}, \alpha_{a_1}), \tau_1(Y_{a_2},
  \alpha_{a_2})\}$ is by construction a canonical base of the type
  $\tp(\tau_1(Y_{z}, \alpha_{z}), \tau_1(Y_{u}, \alpha_{u})/B)$ in the
  full Zariski structure. Since the natural morphism
  $S_0 \to \Aut(k[\eps]/(\eps^2)), s \mapsto \tau_1(Y_s, \alpha_s)$ has
  finite fibres, a canonical base of $\tp(z,u/a_1,a_2$ is
  inter-algebraic with $\{a_1, a_2\}$ in the full Zariski
  structure. Since passing to the reduct can only enlarge a canonical
  base, the canonical base of $\tp(z,u/a_1,a_2)$ is inter-algebraic
  with $\{a_1, a_2\}$. The same argument applies to $\tp(u,v/b_1,b_2)$
  and $\tp(z,v/c_1,c_2)$.

  By Fact~\ref{fieldconf}, the group $G$ is isomorphic to the affine
  group $\G_a(k) \rtimes \G_m(k)$ of an infinite definable field $k$.
\end{proof}

In order to apply the above results we need the following, which is well known model theoretic folklore. We give a short proof specialised to the case where we need it: 
\begin{lemma}
	Let $G$ be a strongly minimal group interpretable in $\Mm$. Then there exists a strongly minimal $Z\subseteq G^2$ that is not a finite boolean combination of cosets of definable subgroups. 
\end{lemma}
\begin{proof}

   To simplify the discussion let us call subsets of $G$ that are
   finite boolean combinations of cosets of $G^n$ (any $n$) affine. By
   strong minimality $G$ is in finite-to-finite correspondence with
   $\Mm$ (this follows, in general, from the fact that $\Mm$ is
   unidimensional. In the present setting $G$ can be assumed to have
   been obtained from Theorem \ref{grconf} using a 1-dimensional group
   configuration, so the existence of a finite-to-finite correspondence
   follows from the statement). It follows that $G$ is not locally
   modular, as the image of any ample family of 1-dimensional subsets
   of $M^2$ under this finite-to-finite correspondence is an ample
   family in $G^2$.
	
   Since $G$ is not locally modular it admits, by \cite[Proposition
   3.21]{ElHaPe} a nearly faithful ample family of generically
   strongly minimal curves $X\to T\subseteq G^2\times T$ of dimension
   3 (i.e. $\dim(T)=3$). Let $G^0$ denote the $\M$-connected component
   of $G$. let $t\in T$ and $x_0\in G^0$ be independent
   $\M$-generics. Let $y_0$ be such that $(x_0,y_0)\in X_t$ and assume
   that $y_0\in gG^0$ for some $g\in G$ (that we can choose
   independent from $(x_0,y_0)$). Then
   $gX_t:=\{(x,y): (x, gy)\in (G^0)^2\cap X_t\}$ is an $\M$-definable
   curve in $(G^0)^2$ and $gX:=\{gX_t: t\in T\}$ is a definable family
   of curves in $(G^0)^2$. Since $G/G^0$ is finite and $X$ is nearly
   faithful the correspondence $X_t\mapsto gX_t$ is finite-to-one and
   on a generic subset of $T$. Therefore $gX$ is readily checked to be
   a 3-dimensional, nearly faithful ample family of curves in
   $(G^0)^2$. Moreover, if $X_t$ is affine (for some $t\in T$) then so
   is $gX_t$. So it will suffice to show that $X$ can be chosen so
   that $gX_t$ is not affine for generic $t\in T$.
	
   If $G^0$ is $\M$-definably isomorphic to either $\G_m$ or to
   an elliptic curve, $\mathcal E$, then for generic $t\in T$
   $gX_t$ is not affine, since there are no definable families of
   subgroups of $\G_m^2$ or of $\mathcal E^2$. Let us elaborate:
   assume towards a contradiction that for generic $t\in T$ the
   curve $X_t$ is affine. So $gX_t$ is also affine for all such
   $t$. In this setting there are finitely many $\0$-definable
   $1$-dimensional subgroups $H_1,\dots H_k$ of $(G^0)^2$ such
   that $gX_t$ coincides, up to a finite set, with a union of
   cosets of the $H_i$. Since $(G^0)^2/H_i$ is $1$-dimensional
   for all $i$, near faithfulness of $gX$ implies that $gX$ is,
   at most, $1$-dimensional, a contradiction.
	
   So we are reduced to the case where $G^0$ is $\M$-definably
   isomorphic to $\G_a$. It is an easy exercise to verify that the
   definable subgroups of $\G_a^2$ are precisely closed subsets cut
   out by linear equations in $\{x^{p^n}\}_{n=0}^\infty$ and
   $\{y^{p^m}\}_{m=0}^\infty$ (for $p=\mathrm{char}(K)$). Thus, if we
   choose $X$ whose projections on $(G^0)^2$ are generically
   unramified (as provided by the combination of Lemma
   \ref{cancel-insep} and Lemma \ref{fr-cancel}) then also $gX$ has
   this property (for a suitable choice of $g$), so if $gX_t$ is
   affine for generic $t$ $gX_t$ contains (up to a finite set) the
   graph of a linear function $x\mapsto ax+b$. As a above, near
   faithfulness of $X$ implies near faithfulness of $gX$, and
   therefore $gX$ can be at most 2-dimensional, again, a
   contradiction.

	
\end{proof}

\begin{rem}
	\begin{enumerate}
		\item It follows from \cite[Theorem 4.1(b)]{hr-pillay} that there is some definable $Z\subseteq G^n$ (some $n$) that is not affine. Reducing $n$ to be 2 requires a little more effort. 
		\item In the above proof it is not hard to see that if we obtain a 2-dimensional nearly faithful family $X\to T$ such that each $X_t$ contains (up to a finite set) a curve of the form $a_tx+b_t$ then $X$ can be used directly to construct a  field configuration. 
	\end{enumerate}
\end{rem}

We can now sum up everything to obtain the main result of this paper: 

\begin{thm}
  \label{main-theorem}
  Let $M$ be an algebraic curve and let
  $X \subset M^2 \times T \subset M^2 \times M^l$ be an ample family of
  curves. Then $\Mm=(M,X)$ interprets a field.
\end{thm}

\begin{proof}
  By Corollary~\ref{reduction}, we may assume that that $M$ is smooth,
  that $k$ is of infinite transcendence degree and that $X$ is a
  nearly faithful family of generically strongly minimal sets. Thus we
  can apply Theorem~\ref{first-group}, allowing us to conclude that
  $\Mm$ interprets a strongly minimal group $G$. By \cite[Theorem
  4.13]{poizat}, $G$ is an algebraic group. The group $G$ is in an
  $\Mm$-definable finite-to-finite correspondence with $M$, so it is a
  one-dimensional algebraic group. Moreover, any $\Mm$-definable ample
  family of curves in $M^2$ maps through this correspondence into an
  ample family of curves in $G^2$ of the same dimension. So $G$ is not
  locally modular. By \cite[Theorem 4.1(b)]{hr-pillay} there is some
  definable $Z\subseteq G$ that is not a finite boolean combination of
  cosets. Let us be a little more precise about this last statement:
  Since $G$ is not locally modular it is not 1-based (or weakly
  normal, in the terminology of \cite{hr-pillay}, see. for example
  \cite[Theorem 8.2.15]{marker}). So by the result of \cite{hr-pillay}
  just referred to one-dimensional, this implies that a generic member
  of a three-dimensional nearly faithful family of curves is not a
  Boolean combination of cosets of algebraic subgroups of $G^2$
  (because any family of cosets of definable subgroups of $G^2$ can be
  -- by dimension considerations -- at most
  two-dimensional). Therefore, we may apply Theorem~\ref{second-group}
  to get the desired conclusion.
\end{proof}

\begin{rem}
  This field is definably isomorphic to $k$ by \cite[Theorem
  4.15]{poizat}.
\end{rem}

\noindent\textit{Acknowledgments.}  The second author thanks Boris
Zilber for his remarks on an early version of the paper, and Maxim
Mornev for many helpful comments. We would like to thank Moshe
Kamensky for comments and suggestions, as well as Ilya Tyomkin and
Martin Bays who, by reading the manuscript, helped us significantly
improve it.

\bibliography{reducts}

\vspace{5ex}

\noindent {\sc Department of Mathematics \\
Ben Gurion University of the Negev \\
P.O.B. 653 \\
Be'er Sheva 84105 \\
Israel\\}
{\tt hassonas@math.bgu.ac.il}
\vspace{2ex}

\noindent {\sc Department of Mathematics\\
  KU Leuven\\
  Celestijnenlaan 200B \\
  B-3001 Leuven (Heverlee)\\
  Belgium\\}
{\tt dmitry.sustretov@kuleuven.be\\}

\end{document}